\documentclass[12pt,a4paper]{book}
\usepackage{amsmath}
\usepackage{amsfonts}
\usepackage{amsthm}
\usepackage{amssymb}
\usepackage{wasysym}
\usepackage{mathtools}
\usepackage{bbm}
\usepackage{tikz-cd}
\usetikzlibrary{babel}

\usepackage{graphics}
\usepackage{fancyhdr}

\swapnumbers

\newtheorem{theorem}{Theorem}[section]

\theoremstyle{definition}

\newtheorem{lemma}[theorem]{Lemma}
\newtheorem{question}[theorem]{Question}
\newtheorem{proposition}[theorem]{Proposition}
\newtheorem{example}[theorem]{Example}
\newtheorem{definition}[theorem]{Definition}
\newtheorem{corollary}[theorem]{Corollary}
\newtheorem{remark}[theorem]{Remark}
\title{Introduction to Contact Dynamics}
\author{Senne Ignoul}

\usepackage{hyperref}

\hypersetup{
    colorlinks=true,
    linkcolor=blue,
    filecolor=blue,      
    urlcolor=blue,
    citecolor = blue,
    }
 
\begin{document}
\pagestyle{fancy}
\renewcommand{\chaptermark}[1]{\markboth{#1}{#1}}
\fancyhead[R]{\thepage}
\fancyhead[LO]{\rightmark}
\fancyhead[LE]{}
\cfoot{}

\thispagestyle{empty}
\begin{titlepage}
\begin{center}
\huge{\bf{\em Introduction to Contact Complete Integrability}}
\normalsize
\ \\
\ \\
\ \\
\ \\
\ \\
Additional lecture notes for the Master course\\ ``Integrable Hamiltonian Systems''
\ \\
\ \\
\ \\
\ \\
\ \\
\resizebox{10cm}{!}{\includegraphics{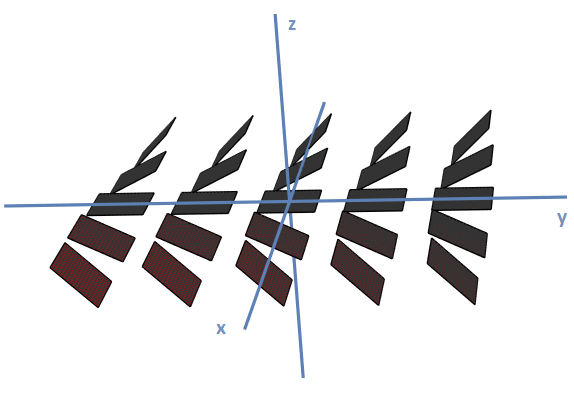}} 
\ \\
\ \\

\ \\
\ \\

\large{Senne Ignoul}\\ 
\large{Universiteit Antwerpen} \\\normalsize
\ \\
\ \\
\ \\
A{\sc cademic year 2023-2024}\\
\ \\
\ \\
\ \\
The author is a predoctoral fellow fundamental research of the Research Foundation - Flanders (FWO),
file number 1179322N
\end{center}
\end{titlepage}
 \thispagestyle{empty}
\tableofcontents
\pagenumbering{roman}
\chapter*{Author's note}
Currently, I am a PhD student at the University of Antwerp. As part of my PhD, I am teaching exercise sessions in a class on Integrable Hamiltonian Systems in our Master Fundamental Mathematics. In this course, we consider Hamiltonian systems in the context of symplectic geometry. Since my own research concerns contact geometry, I decided to also give an introduction in contact integrable systems, which resulted in the current lecture notes.\\
In order for students to be able to understand the literature as soon as possible, my intention is to go to the definition of contact completely integrable systems  as quick as possible. This means that I have not introduced some  notions (such as the contact quotient  or isoenergetic integrability) which are interesting for the study of contact integrable systems, but which are not necessary in order to understand the fundamental notions. In a next edition, I would like to expand these notes by elaborating on these topics.\\
Eventually, prof. dr. Sonja Hohloch (who teaches the abovementioned class) and I want to publish her lecture notes on Integrable Hamiltonian Systems in the symplectic case and my notes on contact complete integrability together as a book (mainly aimed towards graduate students). Hence in these lecture notes, I will often refer to her lecture notes as \textit{the lecture notes on symplectic geometry}. These notes are available on her website:
\begin{center}
\url{https://www.uantwerpen.be/nl/personeel/sonja-hohloch/private-webpage/teaching/antwerpen-this-ye/}.
\end{center}
Finally, I welcome all feedback on these notes. If you have any feedback or if you want to discuss the content of these notes, then please contact me via email \href{mailto:senne.ignoul@uantwerpen.be}{senne.ignoul@uantwerpen.be}.
Finally, I would like to thank the Research Foundation - Flanders (FWO), of which I am a predoctoral fellow fundamental research (1179322N).
\begin{flushright}
Senne Ignoul\\
02/05/2023
\end{flushright}
\chapter{Introduction}
In this course, \textit{Integrable Hamiltonian Systems}, we have discussed a special type of dynamical systems. In particular, we studied completely integrable Hamiltonian systems in the context of symplectic geometry: the symmetries of the system gave us a toolbox to describe the behaviour locally and globally. However, an intrinsic property of symplectic geometry is that it only lives on even dimensional manifolds. This imposes a rather big restriction to the application of this formalism to
concrete problems. Indeed, in many branches of physics (such as thermodynamics), we are interested in solving differential equations on odd dimensional spaces.\\
There exists an odd dimensional sister of symplectic geometry which is called \textit{contact geometry}. In the past, this topic was mostly studied in the context of pure differential geometry and topology. However, in the last few years, there has been an increased interest to use contact geometry to deal with differential equations on odd dimensional spaces. These notes serve as an introduction to these recent developments.\\
The structure of the notes is as follows: in chapter $2$, we will introduce the basic principles of contact geometry (without considering any dynamics). This first part of the chapter is crucial to understand the rest of the notes. In the second part of this chapter, we will see that contact geometry and symplectic geometry are closely related. Indeed, we will discuss the procedure of contactisation and symplectisation.\\
In chapter $3$, we will construct a one-to-one correspondence between smooth functions on a contact manifold and so called contact vector fields (a special class of vector fields which will be introduced in chapter 2). This will be the so-called contact Hamiltonian. Furthermore, we will introduce the Jacobi-bracket, which is different from a Poisson bracket in symplectic geometry, but which has essentially the same dynamical properties (on a special subclass of functions). The algebraic properies of this bracket will be proven in detail, since it shows how one can use this construction in more advanced proofs.\\
Chapter $4$ can be considered as the core of these notes. In this chapter, we will define the notion of a contact completely integrable system, which is our main object of study. However, one needs to know what a (co-)Legendrian foliation is before it is possible to understand this definition. Hence in the first part of chapter $4$, we will introduce these foliations.\\
Since contact dynamics are a relatively new phenomenon, there does not yet exist a standard definition of a contact completely integrable system in the literature. This is why we will give two different definitions, after which we study the relation between them in great detail.\\
While chapter $1$ to $4$ can be considered the contact counterpart of chapter $1$ from the lecture notes on symplectic geometry, chapter $5$ corresponds to chapter $2$ from these lecture notes. Here, we will discuss shortly the semi-local behaviour of a contact completely integrable system near regular points (Arnold-Liouville theorem) and near singular points (Local normal form theorem).\\
In the final chapter (the contact version of chapter $3$ from the lecture notes on symplectic geometry), we will talk about the global aspects of a contact system. In particular, we will discuss the properties of a (contact) moment map for a general Lie group and introduce contact toric $G$-manifolds. We will end with a Delzant-like classification of compact connected contact toric $G$-manifold, which has been proved by Lerman.\\
Due to time limitations (and a lack of knowledge on physics by the author), we will only focus on the mathematical side of the story and ignore the applications in thermodynamics. We refer the more physics-minded reader to the literature, with in particular the work of Alessandro Bravetti (see for example 
\cite{Bravetti}). However, the importance of contact geometry within thermodynamics cannot be underestimated, as illustrated by the following quote of Vladimir Arnold (the same
Arnold from the Arnold-Liouville theorem):
\begin{quote}
Every mathematician knows it is impossible to understand an elementary course in thermodynamics. The reason is that thermodynamics is based - as Gibbs has explicitly proclaimed - on a rather complicated mathematical theory, on contact geometry.
\end{quote}
\chapter{Contact geometry}\label{Chapter Contact geometry}
\pagenumbering{arabic}
In symplectic geometry, the most important object of study is a manifold $W$ together with a symplectic form $\omega$. We have seen that $W$ is necessarily even dimensional. On odd dimensional manifolds, it is possible to study a geometry which is closely related to symplectic geometry, namely contact geometry. In this chapter, we will introduce the basic notions of this geometry, including some additional structures which we will need in later chapters. At the end, we will introduce a method to go from a symplectic manifold to a contact manifold and vice versa. A good reference for a general introduction into contact geometry (without considering dynamics) is \cite{Geiges}, especially the first two chapters.
\section{Basic concepts}
\noindent
In the rest of this section, let $M$ be a $(2n+1)$-dimensional manifold. Instead of a $2$-form, contact geometry deals with a special kind of codimension 1 distribution on $M$ (remember that codimension 1 means that the distribution is 1 dimension lower than the total dimension of $M$, so it is $2n$-dimensional).
\newline

\noindent
\fbox{ \parbox{\textwidth}{
\begin{proposition}
A codimension 1 distribution $\xi$ on $M$ can locally be written as the kernel of a $1$ form on $M$, i.e. for every point $p$ of $M$ there exists an open neighbourhood $U$ of $M$ and a $1$-form $\alpha$ on $M$ such that $\xi|_U = \ker(\alpha|_U)$.
\end{proposition}
}}
\begin{proof}
The proof, which uses some techniques from Riemannian geometry, can be found in \cite{Geiges} (Lemma 1.1.1).
\end{proof}
\noindent
The $1$-form defined above is only defined locally. It would be nice if it can also be defined on the entire manifold, which will be the subject of the next proposition. However, we first have to consider the following new concept.
\newline

\noindent
\fbox{ \parbox{\textwidth}{
\begin{definition}
A codimension 1-distribution $\xi$ on $M$ is called \textit{coorientable} if the quotient $TM/\xi$ (which is a 1-dimensional distribution on $M$) is orientable.
\end{definition}
}}
\newline \hfill \newline
\hfill \newline
\noindent
\fbox{ \parbox{\textwidth}{
\begin{proposition}
Let $\xi$ be a codimension 1 distribution on $M$. There exists a $1$-form $\alpha$ on $M$ such that $\xi = \ker(\alpha)$ if and only if $\xi$ is coorientable.
\end{proposition}}}
\begin{proof}
See Lemma 1.1.1 in \cite{Geiges}.
\end{proof}
\noindent
Recall the definition of the annihilator of $\xi$:
\begin{equation*}
\xi^\circ = \{\beta \in T^*M | \forall u \in \xi: \beta(u) = 0\}. 
\end{equation*}
It can be shown that $\xi^\circ$ is isomorphic to $TM/\xi$. Assume that $\xi$ is coorientable, so there exists a global $1$-form $\alpha$ which defines $\xi$. The coorientation of $\xi$ (which corresponds to an orientation of $\xi^\circ$) will be denoted as $\xi^\circ_+$. This corresponds to the component of $\xi^\circ \setminus \{0\}$ such that $\alpha(M) \subset \xi^\circ_+$ (morally, $\xi^\circ$ is a line in every point and $\alpha$ determines what the positive side of this line is).\newline

\noindent
A well studied class of codimension 1 distributions $\xi$ are the (Frobenius) integrable distributions, i.e. if $X,Y \in \xi$, then $\left[X,Y\right] \in \xi$. If $\xi = \ker(\alpha)$, then it holds that $\xi$ is integrable if and only if $$\alpha\wedge d\alpha = 0.$$ Morally, contact structures will be the exact opposite of these kind of distributions: they will be maximally non-integrable.
\newpage
\noindent
\fbox{ \parbox{\textwidth}{
\begin{definition}
Let $M$ be a $(2n+1)$-dimensional manifold. A differential $1$-form $\alpha$ is called a \textit{contact form} if it holds that
\begin{equation*}
\alpha\wedge (d\alpha)^n \neq 0.
\end{equation*}
A codimension $1$ distribution $\xi$ of $M$ is called a \textit{contact structure} if $\xi$ is (locally) defined by a contact form $\alpha$, so for every $p \in M$, there exists an open neighbourhood $U$ of $p$ such that $\xi|_U = \ker \alpha|_U$.\\
The pair $(M,\xi)$ is called a \textit{contact manifold}. If $\alpha$ is a contact form on $M$, we call the pair $(M,\alpha)$ a \textit{strict contact manifold}.
\end{definition}}}
\noindent
To keep notation at a minimum, we will assume in the continuation of these notes that $\xi$ is always defined by a global contact form $\alpha$ on $M$ (unless stated otherwise). Let us now first consider an example.
\begin{example}
Consider the smooth manifold $\mathbb{R}^{2n+1}$, endowed with Cartesian coordinates $(z,x_1,y_1,\cdots,x_n,y_n)$. Consider the differential $1$-form $\alpha_{st}$ on $\mathbb{R}^{2n+1}$, defined as
\begin{equation*}
\alpha_{st} = dz + \sum_{j=1}^n x_jdy_j.
\end{equation*}
It can be shown by direct calculations (this is mostly a task of bookmarking the indices) that $$\alpha_{st} \wedge (d \alpha_{st})^n = n! dz \wedge
\bigwedge_{j=1}^n  (dx_j \wedge dy_j) \neq 0.$$
Hence $\alpha_{st}$ is a contact form on $\mathbb{R}^{2n+1}$, which is called the \textit{standard contact form} on $\mathbb{R}^{2n+1}$. The corresponding contact structure $\xi = \textrm{Ker}(\alpha_{st})$ is called the \textit{standard contact structure} on $\mathbb{R}^{2n+1}$. Figure $1$ shows the standard contact structure on $\mathbb{R}^3$, which illustrates why this distribution is called maximally non-integrable.
\begin{figure}[h]
\begin{center}
\includegraphics[scale=0.75]{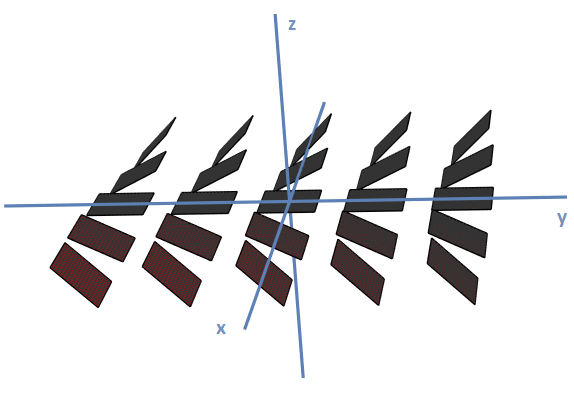}
\caption{The contact structure $\ker(dz + xdy)$ (source \cite{Geiges})}
\end{center}
\end{figure}
\end{example}\noindent
Another way to describe strict contact manifolds is by making use of symplectic vector spaces.
\newline
\noindent
\fbox{ \parbox{\textwidth}{
\begin{proposition}
Let $\xi$ be a codimension $1$ distribution on $M$ (not necessarily contact) that is defined as the kernel of a differential $1$-form $\alpha$. It holds that $\alpha \wedge(d\alpha)^n \neq 0$ if and only if $(d\alpha)^n |_\xi \neq 0$.
\end{proposition}}}
\begin{proof}
This follows by direct calculations.
\end{proof}
\noindent
\fbox{ \parbox{\textwidth}{
\begin{corollary}
$(M,\alpha)$ is a strict contact manifold if and only if $(\xi_m, (d\alpha)_m)$ is a symplectic vector space for every $m \in M$.
\end{corollary}}}
\begin{proof}
This follows directly from the previous proposition, where we have used that the skew-symmetric bilinear form $(d\alpha)_m$ on the vector space $\xi_m$ is non-degenerate if and only if $(d\alpha)_m^n |_{\xi_m} \neq 0$.
\end{proof}\noindent
Associated to every contact form, one can define a special vector field on the contact manifold called the Reeb vector field. This vector field will play a crucial role within the study of contact dynamical systems.
\newline
\noindent
\fbox{ \parbox{\textwidth}{
\begin{definition}
Let $(M,\alpha)$ be a strict contact manifold. The \textit{Reeb vector field} $R_\alpha$ on $M$ associated to $\alpha$ is uniquely defined by
\begin{itemize}
\item[i)] $\iota_{R_\alpha}d \alpha = d \alpha (R_\alpha,\cdot) = 0$,
\item[ii)] $\alpha(R_\alpha) = \mathbbm{1}$
\end{itemize}
The flow of the Reeb vector field is called the \textit{Reeb flow}.
\end{definition}}}
\begin{example}
Consider the smooth manifold $\mathbb{R}^{2n+1}$, endowed with Cartesian coordinates $(z,x_1,y_1,\cdots,x_n,y_n)$ and its standard contact form $\alpha_{st}$. The Reeb vector field $R_{\alpha_{st}}$ of $\alpha_{st}$ is \begin{equation*}
R_{\alpha_{st}} = \frac{\partial}{\partial z}.
\end{equation*}
Indeed, let us check the two conditions:
\begin{itemize}
\item[i)] $d\alpha_{st}(\frac{\partial}{\partial z},\cdot) = d(dz + \sum_{j=1}^n x_jdy_j)(\frac{\partial}{\partial z},\cdot) = \left(\sum_{j=1}^n dx_j \wedge dy_j\right)(\frac{\partial}{\partial z},\cdot) = 0$\\
\item[ii)] $\alpha_{st}\left(\frac{\partial}{\partial z}\right) = dz\left(\frac{\partial}{\partial z}\right) + \sum_{j=1}^n x_jdy_j\left(\frac{\partial}{\partial z}\right) = \mathbbm{1} + \sum_{j=1}^n x_j\cdot 0 = \mathbbm{1}.$
\end{itemize}
\end{example}
\begin{remark}
It is important to notice that the Reeb vector field is an object that heavily depends on the contact form and not only on the contact structure. Indeed, if $\alpha_1$ and $\alpha_2$ are two contact forms which define the same contact structure $\xi$ (so $\xi = \ker(\alpha_1) = \ker(\alpha_2)$), then there needs to exist a function $\lambda: M \rightarrow \mathbb{R}\setminus\{0\}$ such that $\alpha_2 = \lambda \alpha_1$ (they induce the same coorientation if and only if $\lambda$ is positive). However, one can easily check that $R_{\alpha_1} \neq R_{\lambda \alpha_1}$ (provided $\lambda$ is not the constant $1$ function, but then $\alpha_1 = \alpha_2$).
\end{remark}\noindent
Similar to what we have seen in symplectic geometry, all contact manifolds will look locally the same, which is given by the theorem of Darboux. Before we can look at this theorem, we need to introduce a concept which will take on the role of symplectomorphisms from symplectic geometry.
\newline\hfill \newline
\noindent
\fbox{ \parbox{\textwidth}{
\begin{definition}
Let $(M_1,\xi_1)$ and $(M_2,\xi_2)$ be contact manifolds. A diffeomorphism $\varphi: M_1 \rightarrow M_2$ is called a \textit{contactomorphism} if $T\varphi(\xi_1) = \xi_2$. The pairs $(M_1,\xi_1)$ and $(M_2,\xi_2)$ are called contactomorphic.\\
Let $\alpha_1$ and $\alpha_2$ be contact forms which induces the contact structure $\xi_1$ and $\xi_2$, respectively. A diffeomorphism $\varphi: M_1 \rightarrow M_2$ is called a \textit{strict contactomorphism} if $\varphi^*\alpha_2 = \alpha_1$.
\end{definition}}} \newline \hfill \newline \noindent
We can now state the Darboux Theorem for strict contact manifolds.
\newline
\hfill \newline
\noindent
\fbox{ \parbox{\textwidth}{
\begin{theorem}
\textbf{(Darboux) }Let $M$ be a $(2n+1)$-dimensional manifold and let $\alpha$ be a contact form on $M$. For every $p \in M$, there exists an open neighbourhood $U \subset M$ of $p$ and a chart $\varphi: U \rightarrow \mathbb{R}^{2n+1}$ with local coordinates $(z,x_1,y_1,\cdots, x_n,y_n)$ on $U$ such that $\varphi(p) = (0,\cdots, 0)$ and
\begin{equation*}
(\varphi^{-1})^*\alpha = dz + \sum_{j=1}^n x_jdy_j.
\end{equation*}
Hence $(U,\alpha|_{U})$ and $(\mathbb{R}^{2n+1}, dz + \sum_{j=1}^n x_jdy_j)$ are strictly contactomorphic.
\end{theorem}}}
\begin{proof}
See Theorem 2.5.1 in \cite{Geiges}.
\end{proof}\noindent
Dynamics on a manifold are described by a vector field (since differential equations can be described by vector fields). If we want to exploit contact geometry in order to make it easier to describe the dynamics, a first thing we can demand of a vector field is that it preserves the contact structure and/or the contact form. These vector fields will be so important, that we will give them a special name.
\newline \hfill \newline
\noindent
\fbox{ \parbox{\textwidth}{
\begin{definition}
Let $(M,\xi = \ker \alpha)$ be a contact manifold and let $X$ be a vector field on $M$ with (local) flow $\psi_t$ (for all $t \in \mathbb{R}$ such that the flow on $M$ is well defined).\\
We call $X$ a \textit{contact vector field} (or \textit{infinitesimal automorphism of $\xi$}) if $T\psi_t(\xi) = \xi$ for all $t$ such that the flow is well-defined (so the vector field preserves the contact structure).\\
We call $X$ a \textit{strict contact vector field} (or \textit{infinitesimal automorphism of $\alpha$}) if $\psi_t^* \alpha = \alpha$ for all $t$ such that the flow is well-defined (so the vector field preserves the contact form).
\end{definition}}}\newline \hfill \newline \noindent
The following proposition gives us an easy\footnote{Provided you are familiar with methods from global analysis.} way to ``calculate" if a given vector field on a contact manifold is a (strict) contact vector field or not.
\newline
\noindent
\fbox{ \parbox{\textwidth}{
\begin{proposition}\label{contact VF met Lie}
Let $(M,\xi = \ker \alpha)$ be a contact manifold and let $X$ be a vector field on $M$.
\begin{itemize}
\item[i)] $X$ is a contact vector field if and only if there exists a smooth function $\mu: M \rightarrow \mathbb{R}$ such that $\mathcal{L}_X\alpha = \mu \alpha$, where $\mathcal{L}_X$ is the Lie derivative to $X$. This condition does not depend on the choice of $\alpha$.
\item[ii)] X is a strict contact vector field if and only if $\mathcal{L}_{X}\alpha = 0$.
\end{itemize}
\end{proposition}}}
\begin{proof}
The idea of this proof is due to \cite{Geiges}.
\begin{itemize}
\item[i)] First, let us show that the condition $\mathcal{L}_X\alpha = \mu \alpha$ is independent of the choice of the contact form. For any other contact form which defines $\xi$, there exists a nowhere vanishing function $\lambda$ such that this contact form equals $\lambda \alpha$ (see the remark above). We have that
\begin{equation*}
\mathcal{L}_X(\lambda\alpha) = \mathcal{L}_X(\lambda)\alpha + \lambda \mathcal{L}_X(\alpha) = X(\lambda)\alpha + \lambda\mu \alpha = ( X(\lambda) + \lambda\mu)\alpha.
\end{equation*}
Since $\widetilde{\mu} = X(\lambda) + \lambda\mu$ is a function as well, we have proven the independence (the function itself depends on the contact form, but the fact there exists a function doesn't).\\
Now, assume that $T\psi_t(\xi) = \xi$, where $\psi_t$ is the (local) flow of $X$. This means that for every value of $t$, there exists a nowhere vanishing function $\lambda_t$ such that $\psi_t^*\alpha = \lambda_t \alpha$. Consequently,
\begin{equation*}
\mathcal{L}_X\alpha = \frac{d}{dt}\Big|_{t=0}(\psi_t^*\alpha) = \frac{d}{dt}\Big|_{t=0}(\lambda_t \alpha) = \mu \alpha,
\end{equation*}
where $\mu = \frac{d}{dt}\big|_{t=0}\lambda_t$.\\
On the other hand, assume that there exists a function $\mu$ such that $\mathcal{L}_X\alpha = \mu \alpha$. This means that\footnote{We take the derivative in an arbitrary point, not necessarily in $t=0$.}
\begin{equation}\label{Lie derivative 1}
\frac{d}{dt}(\psi_t^*\alpha) = \psi_t^*(\mu \alpha) = \psi_t^*(\mu) \psi_t^*\alpha = (\mu \circ \psi_t)\psi_t^*\alpha,
\end{equation}
where we have used that $\psi_{t+t_0} = \psi_{t_0}\circ \psi_t$, thus via a linear transformation in $t$
\begin{equation*}
\frac{d}{dt}\Big|_{t=t_0}(\psi_t^*\alpha) = \psi^*_{t_0}\left(\frac{d}{dt}\Big|_{t=0}(\psi_t^*\alpha) \right) = \psi^*_{t_0}(\mathcal{L}_X\alpha).
\end{equation*}
Consider equation (\ref{Lie derivative 1}) as a differential equation, we can solve it to find
\begin{equation*}
\psi_t^*\alpha = \lambda_t\alpha
\end{equation*}
where $\lambda_t$ is defined as
\begin{equation*}
\lambda_t = \exp\left(\int_0^t(\mu\circ \psi_s)ds \right).
\end{equation*}
Hence $X$ is a contact vector field.
\item[ii)] $X$ is a strict contact vector field if and only if (the flow of) $X$ preserves the contact form $\alpha$. Using the main property of the Lie derivative, this holds if and only if $\mathcal{L}_X \alpha = 0$.
\end{itemize}
\end{proof}
\begin{example}\label{Reeb VF strict contact VF}
Let $(M,\alpha)$ be a strict contact manifold. It holds by Cartan's magic formula that
\begin{equation*}
\mathcal{L}_{R_\alpha}\alpha = d(\iota_{R_\alpha}\alpha) + \iota_{R_\alpha}d\alpha = d(\mathbbm{1}) + 0 = 0 + 0 = 0.
\end{equation*}
Hence the Reeb vector field associated to $\alpha$ is a strict contact vector field. 
\end{example}
\section{Contactisation and symplectisation}
In the second part of this chapter, we will consider a method to construct a contact manifold of dimension $2n-1$ out of a $2n$-dimensional symplectic manifold, which can be called \textit{contactisation}. Conversely, we will discuss the method of \textit{symplectisation}, which lets us construct a symplectic manifold of dimension $2n$ out of a $(2n-1)$-dimensional contact manifold. At the end, we will see a proposition that relates the Hamiltonian flow in the symplectic setting and the Reeb flow in the contact setting.\\
To construct a contact manifold of dimension $2n-1$ out of a $2n$-dimensional symplectic manifold, there needs to exist a special type of vector field on the symplectic manifold.
\newline \hfill \newline
\noindent
\fbox{ \parbox{\textwidth}{
\begin{definition}
Let $(W,\omega)$ be a symplectic manifold. A vector field $Y$ on $W$ is called a \textit{Liouville vector field} if $\mathcal{L}_Y\omega = \omega$.
\end{definition}}}\\ \hfill \\ \noindent
Notice that since $\omega$ is closed, the condition $\mathcal{L}_Y\omega = \omega$ is, by Cartan's magic formula, equivalent to $d(\iota_Y\omega) = \omega$.\newline

\noindent
We can construct a contact manifold out of a symplectic manifold with a Liouville vector field by considering a special hypersurface (a hypersurface is a codimension 1 submanifold, so a submanifold of 1 dimension lower than the original symplectic manifold).
\newline \hfill \newline
\noindent
\fbox{ \parbox{\textwidth}{
\begin{proposition} \label{contactisation}
\textbf{(Contactisation)}
Let $(W,\omega)$ be a symplectic manifold with a Liouville vector field $Y$. Let $S\subset W$ be a hypersurface of $W$ that is transverse to $Y$ (i.e. $Y$ is nowhere tangent to $S$). Then $\alpha = \iota_Y\omega$ is a contact form on $S$, so $(S,\alpha)$ is a strict contact manifold. The hypersurface $S$ which is transverse to the Liouville vector field is called a \textit{hypersurface of contact type}.
\end{proposition}}}
\begin{proof}
Let $\dim(W) = 2n$ for a certain $n \in \mathbb{N}$. We can calculate
\begin{eqnarray*}
\alpha\wedge (d\alpha)^{n-1} &=& (\iota_Y\omega)\wedge (d(\iota_Y \omega))^{n-1}\\
&=& (\iota_Y\omega)\wedge \omega^{n-1}\\
&=& \frac{1}{n} \iota_Y (\omega^n)
\end{eqnarray*}
Since $(W,\omega)$ is a symplectic manifold, it holds that $\omega^n \neq 0$, so we know that if we restrict to the tangent bundle of a hypersurface transverse to $Y$, it will hold that $\iota_Y(\omega^n) \neq 0$. Consequently, $\alpha \wedge (d\alpha)^{n-1} \neq 0$ on the hypersurface $S$, so $(S,\alpha)$ is a strict contact manifold.
\end{proof}
\noindent
The following proposition describes a method called symplectisation, which can be used to construct a $2n$-dimensional symplectic manifold out of a $(2n-1)$-dimensional contact manifold.
\newline \hfill \newline
\noindent
\fbox{ \parbox{\textwidth}{
\begin{proposition}
\textbf{(Symplectisation)}
Let $(M,\xi = \ker \alpha)$ be a $(2n-1)$-dimensional contact manifold. Consider the manifold $\mathbb{R}\times M$ with an $\mathbb{R}$-coordinate $t$ and let $p: \mathbb{R}\times M \rightarrow M$ be the canonical projection. It holds that $(\mathbb{R}\times M, d(e^tp^*\alpha))$ is a $2n$-dimensional symplectic manifold.
\end{proposition}}}
\begin{proof}
See exercise \ref{oef symplectisatie} from this chapter.
\end{proof}
\noindent
Under certain conditions, the Hamiltonian flow (from symplectic geometry) and the Reeb flow (from contact geometry) will be the same up to reparametrisation.
\\ \hfill
\newline
\noindent
\fbox{ \parbox{\textwidth}{
\begin{proposition}
Let $(W,\omega)$ be a symplectic manifold and consider a non-constant smooth function $H: W \rightarrow \mathbb{R}$. Let $S \subset W$ be a codimension $1$ submanifold of $W$ such that 
\begin{itemize}
\item[i)] $S$ is a hypersurface of contact type with respect to a Liouville vector field $Y$ on $W$ (with contact form $\alpha = \iota_Y\omega$).
\item[ii)] $S$ is a level set of the function $H$.
\end{itemize}
The flow of the Reeb vector field $R_\alpha$ is a reparametrisation of the flow of the Hamiltonian vector field $X^H$ on $S$. 
\end{proposition}}}
\begin{proof}
See Lemma 1.4.10 in \cite{Geiges}.
\end{proof}
\newpage
\section{Exercises}
\begin{enumerate}
\item Let $\alpha$ be a $1$-form on a smooth manifold $M$. Show that $\ker \alpha$ is a Frobenius completely integrable distribution on $M$ if and only if $$\alpha \wedge d\alpha = 0.$$

\item \label{oef T3} Let $n\in \mathbb{N}\setminus \left\{0\right\}$. 
\begin{enumerate}
\item Show that
$$\alpha_n = \cos(nt)d\theta_1 + \sin(nt)d\theta_2$$
is a contact form on $\mathbb{S}^1\times \mathbb{T}^2$, where $(t,\theta_1,\theta_2) \in \mathbb{S}^1 \times \mathbb{T}^2$.
\item Describe the Reeb vector field $R_{\alpha_n}$ of the contact form $\alpha_n$ with respect to the basis $\left\{ \frac{\partial}{\partial t},\frac{\partial}{\partial \theta_1},\frac{\partial}{\partial \theta_2}\right\}$.
\end{enumerate}

\item \label{oef U} Endow $\mathbb{R}^3$ with cylindrical coordinates $(r,\theta,z)$ and consider the set $$U = \left\{(r,\theta,z)| r > 0, 0<\theta < 2\pi \right\}\subset \mathbb{R}^3.$$
\begin{enumerate}
\item Show that $$\alpha = \cos(r)dz + r\sin(r)d\theta$$
defines a contact form on $U$.
\item Describe the Reeb vector field $R_{\alpha}$ of the contact form $\alpha$ with respect to the basis $\left\{ \frac{\partial}{\partial r},\frac{\partial}{\partial \theta},\frac{\partial}{\partial z}\right\}$.

\item Let $X$ be a contact vector field on $(U,\alpha)$ of the form $$X = A\frac{\partial}{\partial r} + B\frac{\partial}{\partial \theta} + C\frac{\partial}{\partial z},$$
where $A,B,C \in \mathbb{R}$. Then $X$ is a strict contact vector field.
\end{enumerate}

\item Is every strict contact manifold orientable? Give a proof or a counterexample.

\item Let $(M,\alpha)$ be a strict contact manifold. Show that the Reeb vector field $R_\alpha$ of $(M,\alpha)$ has no fixed points.

\item Consider $\mathbb{R}^3$ endowed with the Cartesian coordinates $(x,y,z)$. Furthermore, consider the $1$-form
$$\alpha = -y dx - y dy + 2 dz$$
on $\mathbb{R}^3$. Show that the map
$$\varphi: \mathbb{R}^3 \rightarrow \mathbb{R}^3: (x,y,z) \mapsto (-y,x+y,2z)$$
is a strict contactomorphism from $(\mathbb{R}^3,\alpha_{st})$ to $(\mathbb{R}^3,\alpha)$, where $\alpha_{st}$ is the standard contact form on $\mathbb{R}^3$.

\item Consider $\mathbb{R}^3$ endowed with Cartesian coordinates $(x,y,z)$. Furthermore, consider $\alpha_{st} = dz + xdy$, the standard contact form on $\mathbb{R}^3$. 
Decide whether the following vector fields on $\mathbb{R}^3$ are contact vector fields:
\begin{enumerate}
\item $X_1 = \frac{\partial}{\partial x} + x\frac{\partial}{\partial z}$.
\item $X_2 = -xe^y\frac{\partial}{\partial x} + e^y\frac{\partial}{\partial y} + 2\frac{\partial}{\partial z}$.
\item $X_3 = \left(-xy\sin(z) - \cos(z) + \frac{2xy}{y^2+1}\right)\frac{\partial}{\partial x} + \ln\left(y^2+1\right)\frac{\partial}{\partial y} + y\cos(z)\frac{\partial}{\partial z}$.
\end{enumerate}
If the vector field is a contact vector field: decide whether it is a strict contact vector field.

\item \label{oef implicatie Liouville} Let $(M,\omega)$ be a symplectic manifold. Let $Y$ be a vector field on $M$ whose flow $\varphi_t$ (wlog defined for every $t \in \mathbb{R}$) satisfies
$$\forall t \in \mathbb{R}: \varphi_t^*\omega = e^t\omega.$$
Show that $Y$ is a Liouville vector field on $(M,\omega)$.

\item Let $n \in \mathbb{N}\setminus \{0\}$. Consider $\mathbb{R}^{2n}$ with local coordinates $(x_1,y_1,\cdots,x_n,y_n)$ and standard symplectic form $\omega_{st} = \sum_{i=1}^n dy_i \wedge dx_i$.
\begin{enumerate}
\item Show that
\begin{equation*}
Y = \frac{1}{2}\big(\sum_{i=1}^n x_i\frac{\partial}{\partial x_i} + y_i \frac{\partial}{\partial y_i}\big)
\end{equation*}
is a Liouville vector field on $(\mathbb{R}^{2n},\omega_{st})$.
\item Show that $\mathbb{S}^{2n-1} \subset \mathbb{R}^{2n}$ is a hypersurface of contact type with respect to the Liouville vector field $Y$.
\item Use the previous exercises to construct a contact form on $\mathbb{S}^{2n-1}$.
\end{enumerate}

\item Let $M$ be an arbitrary smooth manifold. Consider a Riemannian metric $g$ on $M$. Recall that $g$ induces a fiberwise isomorphism 
\begin{equation*}
g^\flat: TM \rightarrow T^*M, v \mapsto g^\flat(v)
\end{equation*}
defined by $g^\flat(v)(w) = g_m(v,w)$ for $v,w \in T_mM$, and that we denote its inverse by $g^\#: T^*M \mapsto TM$. The metric $g$ on $TM$ induces a metric $g^*$ on $T^*M$, defined as
\begin{equation*}
\forall m \in M, \forall \eta_1,\eta_2 \in T_m^*M: g^*_m(\eta_1,\eta_2) = g_m\left(g^\#(\eta_1)g^\#(\eta_2)\right).
\end{equation*}
We define the \textit{unit cotangent bundle with respect to $g$} as
$$ST^*M = \bigcup_{m\in M} ST_m^*M$$
where we have for every $m \in M$:
$$ST_m^*M = \left\{\eta \in T_m^*M| g_m^*(\eta,\eta) = 1 \right\}.$$
Let $\lambda$ be the Liouville form on $T^*M$. Recall that $(T^*M,d\lambda)$ is a symplectic manifold.\\
Use the technique of contactisation to prove that $(ST^*M,\lambda)$ is a contact manifold.

\item Let $(M,\xi = \ker \alpha)$ be a $(2n-1)$-dimensional contact manifold. Consider the manifold $\mathbb{R}\times M$ with an $\mathbb{R}$-coordinate $t$ and let $p: \mathbb{R}\times M \rightarrow M$ be the canonical projection. 
\begin{enumerate}
\item \label{oef symplectisatie} Show that $$(\mathbb{R}\times M, d(e^tp^*\alpha))$$ is a $2n$-dimensional symplectic manifold.
\item Find a hypersurface of contact type of $$(\mathbb{R}\times M, d(e^tp^*\alpha))$$ such that the contactisation (with the induced contact form) is exactly $(M,\alpha)$.
\end{enumerate}
  
\item Are contactisation and symplectisation inverse constructions, i.e. if you combine both operations, will you end up with the same manifold (up to diffeomorphism!)? Explain your answer.

\end{enumerate}

\newpage
\chapter{Contact Hamiltonians}
In this course, we have always studied special differential equations on symplectic manifolds $(W,\omega)$, namely Hamiltonian equations. These are equations of the form $z' = X^f(z)$, where $f$ is a smooth function on $W$ and $X^f$ is a special vector field on $M$ called the Hamiltonian vector field (defined by the formula $\iota_{X^f}\omega = -df$). In this chapter, we will introduce a special vector field on contact manifolds which will take on the role of Hamiltonian vector field. Furthermore, we will introduce a bracket on functions on contact manifolds which will satisfy the useful dynamical properties of the Poisson bracket in the symplectic case. However, before we can do this, we need to introduce a bit of terminology. We will combine the results from \cite{Geiges}, \cite{Jovanovic} and \cite{Jovanovic Jovanovic}.
\section{Contact Hamiltonians}
\noindent
Consider in this chapter a $(2n+1)$-dimensional strict contact manifold $(M,\alpha)$ with associated contact structure $\xi = \ker \alpha$. Since $\xi$ is a codimension $1$ distribution of $M$ and since $R_\alpha$ is not tangent to $\xi$, we find the following decomposition of the tangent bundle $$TM = \xi \oplus \mathbb{R}R_\alpha.$$ The elements of $\xi$ will be called \textit{horizontal vector field}. Hence if $X$ is a vector field on $M$, we can write $$X = (\iota_X\alpha)R_\alpha + \widetilde{X},$$ for a unique horizontal vector field $\widetilde{X}$ (to verify this, apply $\alpha$ to both sides of the equation). In the rest of these notes, the horizontal part of a vector field $X$ will be denoted as $\widetilde{X}$. Dually, we can decompose the cotangent bundle $T^*M$ by taking the annihilators of both components $$T^*M = \xi^\circ \oplus (\mathbb{R}R_\alpha)^\circ.$$ Notice that, since $\xi$ has codimension $1$, its annihilator has dimension $1$. In particular, we have that $\xi^\circ = \mathbb{R}\alpha$. The elements of $(\mathbb{R}R_\alpha)^\circ$ (hence the $1$-forms $\beta$ such that $\beta(R_\alpha) = 0$) are called \textit{semi-basic forms}. Hence every differential $1$-form $\beta$ on $M$ can be written as $$\beta = \left(\iota_{R_\alpha}\beta\right) \alpha + \widehat{\beta},$$ where $\widehat{\beta}$ is a unique semi-basic form on $M$ (to verify this, evaluate both sides in $R_\alpha$). In the rest of these notes, the semi basic part of a $1$-form $\beta$ will be denoted as $\widehat{\beta}$.\newline

\noindent
We can now consider the map $\alpha^\flat: X \mapsto -\iota_Xd\alpha$ which maps a vector field to a semi-basic form. Notice that, if we restrict $\alpha^\flat$ to horizontal vector fields (i.e. vectors fields in $\xi$), it holds by the non-degeneracy of $d\alpha|_\xi$ that $\alpha^\flat$ is invertible. The inverse map will be denoted by $\alpha^\sharp$.\newline

\noindent
Now, we are able to define contact Hamiltonians.
\newline \hfill
\newline
\noindent
\fbox{ \parbox{\textwidth}{
\begin{theorem} \label{Contact Hamiltonian}
Let $(M,\alpha)$ be a strict contact manifold with $\xi = \ker \alpha$. There exists a one-to-one correspondence between contact vector fields (denoted below as $X$) on $(M,\alpha)$ and smooth functions (denoted below as $H$) on $M$. This correspondence is given by
\begin{itemize}
\item[$\bullet$] $X \longrightarrow H_X = \alpha(X)$
\item[$\bullet$] $H \longrightarrow X_H$, defined by $\alpha(X_H) = H$ and $\iota_{X_H}d\alpha = dH(R_\alpha)\alpha - dH$. 
\end{itemize}
\end{theorem}}}
\begin{proof}
Before we prove the one-to-one correspondence, we will prove that the rules in the second bullet define a unique contact vector field. First of all, if $H$ is a smooth function on $M$, we have that
\begin{equation*}
dH(R_\alpha)\alpha(R_\alpha) - dH(R_\alpha) = dH(R_\alpha) - dH(R_\alpha) = 0,
\end{equation*}
so $dH(R_\alpha)\alpha - dH$ is a semi-basic form. Hence there exists a unique horizontal vector field $\widetilde{X_H} = \alpha^\sharp(dH - dH(R_\alpha)\alpha)$ such that $d\alpha(\widetilde{X_H},\cdot) = dH(R_\alpha)\alpha - dH$. Since $d\alpha(R_\alpha,\cdot) = 0$, we also have for every smooth function $f$ on $M$ that $$d\alpha(fR_\alpha + \widetilde{X_H},\cdot) = dH(R_\alpha)\alpha - dH.$$
Furthermore, we also find, since $\widetilde{X_H}$ is horizontal, that
$$\alpha(fR_\alpha + \widetilde{X_H}) = f \alpha(R_\alpha) + \alpha(\widetilde{X_H}) = f.$$
If we take $f = H$, we have constructed the vector field $HR_\alpha + \widetilde{X_H}$ which satisfies both conditions. By the above considerations, this vector field has to be unique.\\
Using Cartan's magic formula, we can easily prove that $X_H$ is a contact vector field
\begin{equation*}
\mathcal{L}_{X_H}\alpha = d(\iota_{X_H} \alpha) + \iota_{X_H}d\alpha = dH + \left( dH(R_\alpha)\alpha - dH \right) = dH(R_\alpha)\alpha.
\end{equation*}
Consequently, by proposition \ref{contact VF met Lie}, we have that $X_H$ is a contact vector field.\\
Finally, we have to prove that the given correspondence is indeed a one-to-one correspondence. Let $X$ be a contact vector field on $M$. Consider the corresponding smooth function $H_X = \alpha(X)$. We have by Cartan's magic formula that
\begin{equation*}
dH_X + \iota_{X}d\alpha = d(\iota_X\alpha) + \iota_Xd\alpha = \mathcal{L}_X\alpha = \mu \alpha,
\end{equation*}
for a certain smooth function $\mu: M \rightarrow \mathbb{R}$, since $X$ is a contact vector field. If we apply this to $R_\alpha$, we find that:
\begin{equation*}
dH_X(R_\alpha) = (dH_X + \iota_{X}d\alpha)(R_\alpha) = \mu \alpha(R_\alpha) = \mu.
\end{equation*}
Since $\mu = dH_X(R_\alpha)$, it holds that $\iota_{X}d\alpha = dH_X(R_\alpha) \alpha - dH_X$. Since we also know that $\alpha(X) = H_X$, we have proven that $X_{H_X} = X$.\\
Let $H: M \rightarrow \mathbb{R}$ a smooth function. Since by definition, it holds that $\alpha(X_H) = H$, we have that $H_{X_H} = H$, so we are done.
\end{proof}
\noindent
\fbox{ \parbox{\textwidth}{
\begin{definition}
The vector field $X_H$ where $H$ is a smooth function on a contact manifold is called the \textit{contact Hamiltonian vector field} or \textit{contact Hamiltonian}. The corresponding differential equation
\begin{equation*}
z' = X_H(z)
\end{equation*}
is called the \textit{contact Hamiltonian equation}.
\end{definition}}}
\begin{example}
Let $(M,\alpha)$ be a strict contact manifold and consider the constant function
\begin{equation*}
\mathbbm{1}: M \rightarrow \mathbb{R}: m \mapsto 1.
\end{equation*}
One can verify (since $d\mathbbm{1} = 0$) that $X_\mathbbm{1} = R_\alpha$, with $R_\alpha$ the Reeb vector field associated to the contact form $\alpha$.
\end{example}
\noindent
The current definition of the contact Hamiltonian is implicit. However, if we use the isomorphism $\alpha^\sharp$, we are able to make it explicit.
\newline \hfill
\newline
\noindent
\fbox{ \parbox{\textwidth}{
\begin{corollary}
Let $(M,\alpha)$ be a strict contact manifold and let $H$ be a smooth function on $M$. The contact Hamiltonian of $H$ is equal to
\begin{equation*}
X_H = HR_\alpha + \alpha^\sharp (\widehat{dH}).
\end{equation*}
\end{corollary}}}
\begin{proof}
We know from the proof of theorem \ref{Contact Hamiltonian} that $$X_H = HR_\alpha + \alpha^\sharp(dH - dH(R_\alpha)\alpha).$$
We also know that
\begin{equation*}
dH = (\iota_{R_\alpha}dH)\alpha + \widehat{dH} = dH(R_\alpha)\alpha + \widehat{dH},
\end{equation*}
so we conclude that
\begin{equation*}
X_H = HR_\alpha + \alpha^\sharp(\widehat{dH}).
\end{equation*}
\end{proof}
\begin{remark}
If we plug the contact Hamiltonian into $d\alpha$ and then restrict to the contact structure $\xi$, we will get a result which reminds us of the symplectic Hamiltonian vector field (which shouldn't surprise us, since $(\xi_m,(d\alpha)_m)$ is a symplectic vector space for every $m \in M$). Indeed, we get
\begin{equation*}
d\alpha(X_H,\cdot)|_{\xi} = dH(R_\alpha)\alpha|_\xi - dH|_\xi = - dH|_\xi.
\end{equation*}
\end{remark}
\begin{remark}
A notational remark: it should be clear from the context whether  we are talking about the symplectic Hamiltonian or the contact Hamiltonian. However, to be rigorous, we will use the notation $X^f$ for the symplectic Hamiltonian vector field of $f$. On the other hand, we will write $X_f$ for the contact Hamiltonian vector field of $f$.
\end{remark}\noindent
Let us now consider the contact version of the Poisson bracket, namely the Jacobi bracket.
\newline \hfill
\newline
\noindent
\fbox{ \parbox{\textwidth}{
\begin{definition}
Let $(M,\alpha)$ be a strict contact manifold. We define the \textit{Jacobi-bracket} $[\cdot,\cdot]$ as
\begin{equation*}
[\cdot,\cdot]: \mathcal{C}^{\infty}(M,\mathbb{R})\times \mathcal{C}^{\infty}(M,\mathbb{R}) \rightarrow \mathcal{C}^{\infty}(M,\mathbb{R}): (f,g) \mapsto [f,g] = \alpha([X_f,X_g]),
\end{equation*}
where $[X_f,X_g]$ is the standard Lie bracket of the vector fields $X_f$ and $X_g$.
\end{definition}}}\newline \hfill \newline \noindent 
Notice that, since $\mathcal{L}_{[X_f,X_g]} = \mathcal{L}_{X_f}\circ \mathcal{L}_{X_g} - \mathcal{L}_{X_g}\circ \mathcal{L}_{X_f}$, we have by proposition \ref{contact VF met Lie} that $[X_f,X_g]$ is again a contact vector field on $M$. Hence it is natural to use the contact Hamiltonian correspondence to define a smooth function in this way.\newline

\noindent
The Jacobi bracket satisfies some algebraic properties which have a lot of practical uses. As stated earlier, it is a contact version of the Poisson bracket, but there are some important differences as well. In what follows, the proofs are a bit algebraic in nature, but they will be given nonetheless for the sake of completeness.
\newline \hfill \newline
\noindent
\fbox{ \parbox{\textwidth}{
\begin{proposition}
Let $(M,\alpha)$ be a strict contact manifold. Let $[\cdot,\cdot]$ be the Jacobi-bracket on $\mathcal{C}^{\infty}(M,\mathbb{R})$. It holds that $(\mathcal{C}^{\infty}(M,\mathbb{R}),[\cdot,\cdot])$ is a Lie algebra.
\end{proposition}}}
\begin{proof}
Let $f,g,h \in \mathcal{C}^{\infty}(M,\mathbb{R})$ and let $\lambda\in \mathbb{R}$. First of all, it holds that $[\lambda f + g, h] = \alpha([X_{\lambda f + g},X_h])$. It is easy to see from the correspondence described in proposition \ref{Contact Hamiltonian} that $X_{\lambda f + g} = \lambda X_f + X_g$. Hence we have that
\begin{equation*}
[\lambda f + g, h] = \alpha([\lambda X_{f} + X_{g},X_h]) = \lambda \alpha([X_f,X_h]) + \alpha([X_g,X_h]) = \lambda [f,h] + [g,h].
\end{equation*}
Now, let $f,g \in \mathcal{C}^{\infty}(M,\mathbb{R})$. It holds that
\begin{equation*}
[f,g] = \alpha([X_f,X_g]) = -\alpha([X_g,X_f]) = -[g,f].
\end{equation*}
Finally, in order to prove the Jacobi-identity, let $f,g,h \in \mathcal{C}^{\infty}(M,\mathbb{R})$. We first notice that $\alpha(X_{[g,h]}) = [g,h] = \alpha([X_g,X_h])$ and since both $X_{[g,h]}$ and $[X_g,X_h]$ are contact vector fields, it follows that $X_{[g,h]} = [X_g,X_h]$. Hence
\begin{equation*}
[f,[g,h]] = \alpha([X_f,X_{[g,h]}]) = \alpha([X_f,[X_g,X_h]]).
\end{equation*}
By cyclic permutation of $f,g$ and $h$, this gives us that
\begin{eqnarray*}
0 &=& \alpha([X_f,[X_g,X_h]]+[X_g,[X_h,X_f]] + [X_h,[X_f,X_g]])\\
&=& \alpha([X_f,[X_g,X_h]]) + \alpha([X_g,[X_h,X_f]]) + \alpha([X_h,[X_f,X_g]])\\
&=& [f,[g,h]] + [g,[h,f]] + [h,[f,g]],
\end{eqnarray*}
where we have used the Jacobi-identity of the Lie bracket of vector fields in the first equality. This proves the Jacobi-identity for the Jacobi-bracket and thus, $(\mathcal{C}^{\infty}(M,\mathbb{R}),[\cdot,\cdot])$ is a Lie algebra. 
\end{proof}\noindent
We have already seen that the set of contact vector fields on a contact manifold is closed under the Lie bracket of vector fields. Denote
\begin{equation*}
\mathcal{X}_{con}(M,\xi) = \{X | X \textrm{ is a contact vector field on }(M,\xi)\}
\end{equation*}
for a contact manifold $(M,\xi)$. We have that $(\mathcal{X}_{con}(M,\xi),[\cdot,\cdot])$ with $[\cdot,\cdot]$ the Lie bracket of vector fields is a Lie algebra.
\newline \hfill
\newline
\noindent
\fbox{ \parbox{\textwidth}{
\begin{proposition}
Let $(M,\alpha)$ be a strict contact manifold. The map
\begin{equation*}
\Phi: (\mathcal{C}^{\infty}(M,\mathbb{R}),[\cdot,\cdot]) \rightarrow (\mathcal{X}_{con}(M,\ker \alpha),[\cdot,\cdot]): f \mapsto X_f
\end{equation*}
is a Lie algebra isomorphism.
\end{proposition}}}
\begin{proof}
We already know that $\Phi$ is an isomorphism of $\mathbb{R}$-vector spaces (it is clear that $\Phi$ is $\mathbb{R}$-linear and it is bijective because the correspondence given in proposition \ref{Contact Hamiltonian} is a one-to-one correspondence). Hence the only thing left to prove is that $\Phi([f,g]) = [\Phi(f),\Phi(g)]$ for every $f,g \in \mathcal{C}^{\infty}(M,\mathbb{R})$, or equivalently that $X_{[f,g]} = [X_f,X_g]$. Recall that both $X_{[f,g]}$ and $[X_f,X_g]$ are contact vector fields such that $\alpha(X_{[f,g]}) = [f,g]= \alpha([X_f,X_g])$. Since the correspondence in proposition \ref{Contact Hamiltonian} is a one-to-one correspondence, we have that $X_{[f,g]} = [X_f,X_g]$, which we had to prove.
\end{proof}\noindent
Although the Jacobi-bracket is a Lie bracket on the smooth functions of a contact manifold, it is not a Poisson bracket, since it will not satisfy the Leibniz rule. Before we can prove this, we will give two characterizations of the Jacobi bracket which will be used a lot in applications.
\newline \hfill
\newline
\noindent
\fbox{ \parbox{\textwidth}{
\begin{proposition} \label{Jacobi karakterisatie 1}
Let $(M,\alpha)$ be a strict contact manifold and consider $f,g \in \mathcal{C}^{\infty}(M,\mathbb{R})$. It holds that
\begin{equation*}
[f,g] = d\alpha(X_f,X_g) + f R_\alpha(g) - g R_\alpha(f).
\end{equation*}
\end{proposition}}}
\begin{proof}
We know from differential geometry/global analysis that
\begin{equation*}
d\alpha(X_f,X_g) = (\mathcal{L}_{X_f}\alpha)(X_g) - (\mathcal{L}_{X_g}\alpha)(X_f) + \alpha([X_f,X_g]).
\end{equation*}
Since $[f,g] = \alpha([X_f,X_g])$ and since $\mathcal{L}_{X_h}\alpha = dh(R_\alpha)\alpha$ for every smooth function $h$ (we have seen this in the proof of theorem \ref{Contact Hamiltonian}), it holds that
\begin{eqnarray*}
[f,g] &=& d\alpha(X_f,X_g) + (\mathcal{L}_{X_g}\alpha)(X_f) - (\mathcal{L}_{X_f}\alpha)(X_g)\\
&=&  d\alpha(X_f,X_g) + (dg(R_\alpha)\alpha)(X_f) - (df(R_\alpha)\alpha)(X_g)\\
&=&d\alpha(X_f,X_g) + R_\alpha(g) \alpha(X_f) - R_\alpha(f) \alpha(X_g)\\
&=&d\alpha(X_f,X_g) + f R_\alpha(g) - g R_\alpha(f),
\end{eqnarray*}
which we wanted to prove.
\end{proof}\noindent
Another realization is given by the following proposition.
\newline \hfill
\newline
\noindent
\fbox{ \parbox{\textwidth}{
\begin{proposition}\label{Jacobi karakterisatie 2}
Let $(M,\alpha)$ be a strict contact manifold and consider $f,g \in \mathcal{C}^{\infty}(M,\mathbb{R})$. It holds that
\begin{equation*}
[f,g] = X_f(g)-gR_\alpha(f)
\end{equation*}
\end{proposition}}}
\begin{proof}
Recall from differential geometry/global analysis that
\begin{equation*}
\iota_{[X_f,X_g]} = \mathcal{L}_{X_f}\iota_{X_g} - \iota_{X_g}\mathcal{L}_{X_f}.
\end{equation*}
Hence we find that
\begin{eqnarray*}
[f,g] &=& \iota_{[X_f,X_g]}\alpha\\
&=& \mathcal{L}_{X_f}\iota_{X_g}\alpha - \iota_{X_g}\mathcal{L}_{X_f}\alpha\\
&=& \mathcal{L}_{X_f}g - \iota_{X_g}(df(R_\alpha)\alpha)\\
&=& X_f(g) - df(R_\alpha) \alpha(X_g)\\
&=& X_f(g) - gR_\alpha(f),
\end{eqnarray*}
so we are done.
\end{proof}\noindent
We can now finally prove that the Jacobi bracket does not satisfy the Leibniz rule.
\newline \hfill
\newline
\noindent
\fbox{ \parbox{\textwidth}{
\begin{proposition}
Let $(M,\alpha)$ be a strict contact manifold and let $f,f'$ and $g$ be smooth functions on $M$. It holds that
\begin{equation*}
[ff',g] = f[f',g] + f'[f,g] - ff'[\mathbbm{1},g].
\end{equation*}
\end{proposition}}}
\begin{proof}
First of all, we will prove that
\begin{equation}\label{Hamiltoniaan ff'}
X_{ff'} = fX_{f'} + f'X_f - ff'R_\alpha.
\end{equation}
Verifying the first condition of a contact Hamiltonian, we find that
\begin{eqnarray*}
\alpha(fX_{f'} + f'X_f - ff'R_\alpha) &=& f\alpha(X_{f'}) + f'\alpha(X_f) - ff'\alpha(R_\alpha)\\ 
&=& ff' + f'f - ff'\\ 
&=& ff'.
\end{eqnarray*}
On the other hand, we will find for the second condition that
\begin{eqnarray*}
d\alpha(fX_{f'} + f'X_f - ff'R_\alpha,\cdot) &=& d\alpha(fX_{f'} + f'X_f,\cdot)\\
&=& f d\alpha(X_{f'},\cdot) + f'd\alpha(X_f,\cdot)\\
&=& f\left(df'(R_\alpha)\alpha - df' \right) + f'\left(df(R_\alpha)\alpha - df \right)\\
&=& fR_\alpha(f')\alpha - fdf' + f'R_\alpha(f)\alpha - f'df\\
&=& R_\alpha(ff')\alpha - d(ff'),
\end{eqnarray*}
where we have used the Leibniz rule of vector fields and of the exterior derivative applied on $0$-forms (i.e. functions). Hence we have proven equation (\ref{Hamiltoniaan ff'}).\newline

\noindent
Now, using proposition \ref{Jacobi karakterisatie 2}, we find that
\begin{eqnarray*}
[ff',g] &=& X_{ff'}(g) - gR_\alpha(ff')\\
&=& fX_{f'}(g) + f'X_f(g) - ff'R_\alpha(g) - gR_\alpha(ff')\\
&=& fX_{f'}(g) + f'X_f(g) - ff'R_\alpha(g) - g\left(fR_{\alpha}(f') + f'R_\alpha(f)\right)\\
&=& f\left(X_{f'}(g) -g R_\alpha(f')\right) + f'\left(X_f(g) - gR_\alpha(f)\right) - ff'\left(X_\mathbbm{1}(g) - gR_\alpha(\mathbbm{1})\right)\\
&=& f[f',g] + f'[f,g] - ff'[\mathbbm{1},g],
\end{eqnarray*}
where we have used that $X_\mathbbm{1} = R_\alpha$ and that $R_\alpha(\mathbbm{1}) = 0$.
\end{proof}\noindent
Hence in general, the Jacobi bracket is not a Poisson bracket. However, our main interest in Poisson brackets in this course are the dynamical properties (in particular concerning integrals) that are associated to them. Now, it turns out that if we only consider integrals of the Reeb vector field, the Jacobi bracket will satisfy the useful dynamical properties!\newline

\noindent
Let $(M,\alpha)$ be a strict contact manifold. We will denote the space of integrals of the Reeb vector field $R_\alpha$ as
\begin{equation*}
\mathcal{C}^{\infty}_\alpha(M,\mathbb{R}) = \{f \in \mathcal{C}^{\infty}(M,\mathbb{R}) |\textrm{  } [\mathbbm{1},f] = R_\alpha(f) = 0\}.
\end{equation*}
\newline \hfill
\newline
\noindent
\fbox{ \parbox{\textwidth}{
\begin{proposition}\label{Jacobi bracket als Poisson bracket}
Let $(M,\alpha)$ be a strict contact manifold and consider $f \in \mathcal{C}^{\infty}_\alpha(M,\mathbb{R})$.
\begin{itemize}
\item[a)] If $g \in \mathcal{C}^{\infty}_\alpha(M,\mathbb{R})$, it holds that the following are equivalent
\begin{itemize}
\item[i)] $g$ is an integral of $f$.
\item[ii)] $f$ is an integral of $g$.
\item[iii)] $[f,g] = 0$.
\end{itemize}
\item[b)] For every $g_1,g_2 \in \mathcal{C}^{\infty}_\alpha(M,\mathbb{R})$, it holds that if $g_1$ and $g_2$ are integrals of $f$, so is $[g_1,g_2]$.
\end{itemize}
\end{proposition}}}
\begin{proof}
See exercise \ref{oef jacobi} from this chapter.
\end{proof}
\newpage
\section{Exercises}

\begin{enumerate}
\item Given one of the following strict contact manifolds $(M,\alpha)$ and vector fields $X$ on $M$, determine the unique horizontal vector field $\widetilde{X}$ associated to $X$ through the decomposition $TM = \ker \alpha \oplus \mathbb{R}R_\alpha$:
\begin{enumerate}
\item $(\mathbb{R}^3, dz + xdy)$ with the vector field
$$X_1 = \sinh(x^2 + y) \frac{\partial}{\partial x} + 42 \frac{\partial}{\partial y} - (z + xy^2)\frac{\partial}{\partial z}. $$
\item $\left(\mathbb{S}^1 \times \mathbb{T}^2, \alpha_n\right)$ (see exercise \ref{oef T3} from chapter \ref{Chapter Contact geometry}) with the vector field $$X_2 = e^{\cos \theta_1} \frac{\partial}{\partial t} - \ln\left(\sin^2(nt) + 1\right)\frac{\partial}{\partial \theta_2}.$$
\item $(U,\alpha)$ (see exercise \ref{oef U} from chapter \ref{Chapter Contact geometry}) with the vector field
$$X_3 = z \frac{\partial}{\partial r} + e^z \cos^2(r)\frac{\partial}{\partial \theta} + e^zr^2\frac{\partial}{\partial z}.$$
\end{enumerate}

\item Given one of the following strict contact manifolds $(M,\alpha)$ and $1$-forms $\beta$ on $M$, determine the unique semi-basic form $\widehat{\beta}$ associated to $\beta$ through the decomposition $T^*M = \left(\ker \alpha\right)^\circ \oplus \left(\mathbb{R}R_\alpha\right)^\circ$:
\begin{enumerate}
\item $(\mathbb{R}^3, dz + xdy)$ with the $1$-form
$$\beta_1 =  (x^2+y)dx + \frac{1}{z^2 + 1} dy + xyzdz.$$
\item $\left(\mathbb{S}^1 \times \mathbb{T}^2, \alpha_n\right)$ (see exercise \ref{oef T3} from chapter \ref{Chapter Contact geometry}) with the $1$-form 
$$\beta_2 =\cos^3(\theta_2) dt + e^{\cos(nt)}\cos(nt)d\theta_1 + e^{\cos(nt)}\sin(nt)d\theta_2.$$
\item $(U,\alpha)$ (see exercise \ref{oef U} from chapter \ref{Chapter Contact geometry}) with the $1$-form
$$\beta_3 = zdr +\left(r+\cos(r)\sin(r)\right) \left(r\sin(r)d\theta + \cos(r)dz\right).$$
\end{enumerate}

\item Let $H$ be a smooth function on a strict contact manifold $(M,\alpha)$. Is the contact Hamiltonian vector field $X_H$ a strict contact vector field? If yes, give a proof. If not, can you find a condition on $H$ such that $X_H$ is a strict contact vector field?

%
%

\item \label{oef jacobi} In this exercise, we will prove the special \textit{dynamical} properties of the Jacobi bracket. Consider a strict contact manifold $(M,\alpha)$ with Reeb vector field $R_\alpha$.\\
\hfill  \\
Recall that we say that a smooth function $h$ on $M$ is an integral of the (not necessarily contact) vector field $X$ if $X(h) = 0$. Furthermore, if $h$ and $k$ are smooth functions on $M$, then $h$ is an integral of $k$ if $h$ is an integral of the contact Hamiltonian $X_k$.
\begin{enumerate}
\item Show that for a smooth function $f$ on $M$, it holds that $f$ is an integral of $R_\alpha$ if and only if $\left[1,f\right] = 0$.
\item Let $f,g \in \mathcal{C}^\infty_\alpha(M,\mathbb{R})$. Show that the following are equivalent:
\begin{itemize}
\item[i)] $f$ is an integral of $g$.
\item[ii)] $g$ is an integral of $f$.
\item[iii)] $\left[f,g\right] = 0$.
\end{itemize}
\item Let $f,g_1,g_2 \in \mathcal{C}^\infty_\alpha(M,\mathbb{R})$. Show that if $g_1$ and $g_2$ are integrals of $f$, then so is $\left[g_1,g_2\right]$.
\end{enumerate}

\item Let $(M,\alpha)$ be a strict contact manifold and consider the associated Jacobi bracket $[\cdot,\cdot]$. Does there exist a non-trivial class of functions on which the Jacobi bracket is a Poisson bracket?

\end{enumerate}

\newpage
\chapter{Contact completely integrable systems}
In the literature, there exist many notions of \textit{contact completely integrable systems}, which are not all completely equivalent. In these notes, we will take a look at two different formalisms.\\
The first one was introduced by \cite{Khesin Tabachnikov} and is very geometric in nature. Between the two formalisms, the intuition why this definition is useful should be the clearest. However, there are two big disadvantages: first of all, it uses a very difficult concept from advanced differential geometry (but we will be able to work around this). Secondly, it will not be very easy to use this definition on concrete examples.\\
The second formalism, which was worked out by \cite{Jovanovic Jovanovic} and \cite{Jovanovic} (but \cite{Khesin Tabachnikov} also hints at it), will look a lot more like the concept of completely integrable systems from symplectic geometry. This formalism will be more useful in applications. However, a disadvantage is that you can only use it on strict contact manifolds.\\
Before we introduce both formalisms, let us take a look at some new geometric concepts related to contact geometry.
\section{Legendrian and co-Legendrian foliations}
In symplectic geometry, we have considered special classes of submanifolds of a symplectic manifold (isotropic and Lagrangian submanifolds). In contact geometry, we have similar special submanifolds.
\newline \hfill
\newline
\noindent
\fbox{ \parbox{\textwidth}{
\begin{definition}
Let $(M,\xi)$ be a contact manifold. A submanifold $L$ of $M$ will be called \textit{isotropic} if for every $p \in L$, it holds that $T_pL \subseteq \xi_p$, i.e. the distribution $TL$ is a part of the distribution $\xi|_L$.
\end{definition}}}\newline \hfill \newline \noindent
This definition looks completely different from the notion of isotropic submanifold from symplectic geometry. However, both of them will be related as follows.
\newline \hfill
\newline
\noindent
\fbox{ \parbox{\textwidth}{
\begin{proposition}
Let $(M,\xi)$ be a contact manifold and let $L$ be an isotropic submanifold of $M$. Then for every $p \in L$, we have that $T_pL$ is an isotropic subspace of the symplectic vector space $(\xi_p,d\alpha|_{\xi_p})$ (where $\alpha$ is a contact form that defines $\xi$ locally around $p$).
\end{proposition}}}
\begin{proof}
Consider the natural inclusion $i: L\hookrightarrow M$. Let $p \in L$ and let $\alpha$ be a contact form that defines $\xi$ locally around $p$. Since $L$ is an isotropic submanifold, it holds that $\alpha_p(u) = 0$ for every $u \in T_pL$. We can write this condition as $(i^*\alpha)_p = 0$. Hence we have that
\begin{equation*}
i^*d\alpha_p = d(i^*\alpha)_p = d(0) = 0.
\end{equation*}
This means exactly that $T_pL$ is an isotropic subspace of the symplectic vector space $(\xi_p,d\alpha|_{\xi_p})$, namely for every $u,v \in T_pL$, it holds that $d\alpha_p(u,v) = 0$.
\end{proof}\noindent
This proposition has a strong consequence: the dimension of isotropic submanifolds of a contact manifold is bounded.
\newline \hfill
\newline
\noindent
\fbox{ \parbox{\textwidth}{
\begin{corollary}
Let $(M,\xi)$ be a $(2n+1)$-dimensional contact manifold and let $L$ be an isotropic submanifold of $M$. Then $\dim(L) \leq n$.
\end{corollary}}}
\begin{proof}
Let $p \in L$ and let $\alpha$ be a contact form that defines $\xi$ locally around $p$. We know by the previous proposition that $T_pL$ is an isotropic subspace of the $2n$-dimensional symplectic vector space $(\xi_p,\omega_p = d\alpha|_{\xi_p})$. If $(T_pL)^\omega$ is the symplectic complement of $T_pL$ in $\xi_p$, we have that $T_pL \subseteq (T_pL)^\omega$ (this is another way to write down the isotropy) and
\begin{equation*}
\dim(T_pL) + \dim((T_pL)^\omega) = \dim(\xi_p) = 2n.
\end{equation*}
We conclude that $\dim(T_pL)\leq n$, so $\dim(L) \leq n$.
\end{proof}
\noindent
\fbox{ \parbox{\textwidth}{
\begin{definition}
Let $(M,\xi)$ be a $(2n+1)$-dimensional contact manifold. A submanifold $L$ of $M$ is called \textit{Legendrian} if $L$ is an isotropic submanifold of maximal dimension, so $\dim(L) = n$.
\end{definition}}}\newline \noindent \newline \noindent
In what follows, we are interested in foliations and not so much in submanifolds. Recall that if $M$ is a manifold, a (regular) $k$-dimensional \textit{foliation} $\mathcal{F}$ on $M$ is a collection of disjoint, non-empty, connected, immersed $k$-dimensional submanifolds of $M$ such that the union of all these submanifolds is $M$ (there is an additional condition on the possible charts, but we will refer the reader to \cite{Lee} for a more in depth discussion of  foliations). These submanifolds are called the \textit{leaves} of the foliation. If not all the submanifolds are $k$-dimensional, one refers to $\mathcal{F}$ as a singular $k$-dimensional foliation. Unless stated otherwise, we will always assume that we are working with regular foliations without explicitly mentioning it.
\begin{example}
Consider the smooth manifold $\mathbb{R}^2$. An example of a $1$-dimensional foliation on $\mathbb{R}^2$ is the collection of submanifolds of the form
\begin{equation*}
\forall c \in \mathbb{R}: \{(x,y)\in \mathbb{R}^2| y=c\},
\end{equation*}
which is a collection of horizontal lines. On the other hand, an example of a singular $1$-dimensional foliation on $\mathbb{R}^2$ is the collection of submanifolds of the form (where we work in polar coordinates)
\begin{equation*}
\forall a \in \mathbb{R}^+\cup \{0\}: \{re^{i\theta}| r = a, \theta \in \left[0,2\pi\right[\},
\end{equation*}
which is a collection of concentric circles with radius $a$. If $a = 0$, then we have just the origin, which is a $0$-dimensional submanifold of $\mathbb{R}^2$.
\end{example}\noindent
In order to define contact completely integrable systems, we need to introduce following special kinds of foliations.
\newline \hfill
\newline
\noindent
\fbox{ \parbox{\textwidth}{
\begin{definition}
Let $(M,\xi)$ be a $(2n+1)$-dimensional contact manifold. An $(n+1)$-dimensional foliation $\mathcal{F}$ on $M$ is called a \textit{co-Legendrian} foliation if
\begin{itemize}
\item[i)]$\mathcal{F}$ is transverse to $\xi$, i.e. $\xi$ is not tangent to $\mathcal{F}$ (or more precisely, not tangent to any leaf of $\mathcal{F}$).
\item[ii)]$T\mathcal{F} \cap \xi$ is an integrable distribution (or more precisely, $TF\cap \xi$ is an integrable distribution for every leaf $F$ of $\mathcal{F}$) , i.e. for every $x \in M$, there exists a unique maximal integral manifold $N$ of $T\mathcal{F} \cap \xi$ through $x$ (so for every $y \in N$, it holds that $T_yN = (T\mathcal{F} \cap \xi)|_y)$.
\end{itemize}
\end{definition}}}
\newline \hfill
\newline
\noindent
\fbox{ \parbox{\textwidth}{
\begin{definition}
Let $(M,\xi)$ be a $(2n+1)$-dimensional contact manifold. An $n$-dimensional foliation $\mathcal{G}$ on $M$ is called a \textit{Legendrian} foliation if every leaf of $\mathcal{G}$ is a Legendrian submanifold of $(M,\xi)$.
\end{definition}}}\newline \hfill \newline \noindent
In the rest of this section, we will study some important properties of co-Legendrian foliations on a contact manifold. First of all, we will see that one can associate a Legendrian foliation to a co-Legendrian foliation. Conversely, given the flow of a certain contact vector field, one can construct a co-Legendrian foliation from a Legendrian foliation.\newline

\noindent
Let $\mathcal{F}$ be a co-Legendrian foliation on the $(2n+1)$-dimensional contact manifold $(M,\xi)$. Since $T\mathcal{F} \cap \xi$ is integrable, there exist a foliation $\mathcal{G}$ of the distribution $T\mathcal{F} \cap \xi$ on $M$. Since $T\mathcal{G} = T\mathcal{F} \cap \xi  \subset \xi$, we have that $\mathcal{G}$ is an isotropic foliation (i.e. the leaves are isotropic submanifolds). Furthermore, because $\mathcal{F}$ is transverse to $\xi$, we know that $\dim(T\mathcal{F} \cap \xi) = n$ (in other words, $1$ dimension less than the dimension of $T\mathcal{F}$) and thus we also have that $\dim(\mathcal{G}) = n$. We conclude that $\mathcal{G}$ is a Legendrian foliation (i.e. the leaves are Legendrian submanifolds), which we call the \textit{Legendrian foliation associated to the co-Legendrian foliation} $\mathcal{F}$. Notice that $\mathcal{G}$ is a codimension $1$ foliation of every leaf of $\mathcal{F}$.\newline

\noindent
Conversely, let $(M,\xi)$ be a contact manifold, $\mathcal{G}$ a Legendrian foliation on $(M,\xi)$ and $X$ a contact vector field on $M$ which is transverse to $\xi$ (if we have that $\xi = \ker(\alpha)$, then the Reeb vector field of $\alpha$ is an example of such a vector field). Furthermore, assume that $X$ preserves $\mathcal{G}$ leafwise (so the flow of $X$ sends two points of the same of leaf of $\mathcal{G}$ to the same of leaf of $\mathcal{G}$). If $x \in M$ lies in the leaf $G_1$ of $\mathcal{G}$, then (since $X$ is transverse to $\xi$) the flow of $X$ will send $x$ to a different leaf $G_2$ of $\mathcal{G}$ (if the time of the flow is small enough). The orbits of the leaves of $\mathcal{G}$ under the action induced by $X$ will form a new foliation $\mathcal{F}$ which is $1$ dimension higher then $\mathcal{G}$. Since $X$ is transverse to $\xi$ and since $T\mathcal{F}\cap \xi = T\mathcal{G}$, it is clear that $\mathcal{F}$ is a co-Legendrian foliation on $(M,\xi)$. In fact, lemma $2.2$ from \cite{Khesin Tabachnikov} gives us that every co-Legendrian foliation looks locally like this construction.
\newline \hfill
\newline
\noindent
\fbox{ \parbox{\textwidth}{
\begin{proposition}\label{model co legendrian}
If $\mathcal{F}$ is a co-Legendrian foliation on a contact manifold $(M,\xi)$, then $\mathcal{F}$ is locally contactomorphic to the co-Legendrian foliation as constructed above.
\end{proposition}}}
\newpage
\begin{proof}
We will follow the proof by \cite{Khesin Tabachnikov} very closely. First of all, we define\footnote{This contact manifold is one of the first contact manifolds which has been described and is in fact the reason why it is called a \textit{contact} manifold.} an important contact manifold associated to a general smooth manifold $N$.\\
We define a contact element of $N$ (which has so far nothing to do with contact geometry) as a hyperplane in a tangent space of $N$. The space of contact elements $C(N)$ consists of pairs $(n,V)$ where $n\in N$ and where $V \subset T_nN$ is a contact element of $N$. Notice that $C(N)$ can be seen as the projectivised cotangent bundle $\mathbb{P}T^*N$, which is defined as:
\begin{equation*}
\mathbb{P}T^*N = \{(n,\beta) | n \in N, \beta \in \mathbb{P}T_n^*N \},
\end{equation*}
where we have that $\mathbb{P}T_n^*N$ is the projectivization of the vector space $T_n^*N$ for every $n \in N$. Indeed, a hyperplane in $T_nN$ is defined as the kernel of an element of $T_n^*N$, which is defined uniquely up to multiplication by a non-zero scalar. Hence the identification is given by associating a hyperplane in a tangent space with a vector (up to rescaling) in the cotangent space which becomes $0$ on the hyperplane.\\
Now, consider the canonical projection $p: \mathbb{P}T^*N \rightarrow N$ which sends each contact element to its foot point. For every $v \in \mathbb{P}T^*N$, we define the hyperplane $\xi_v$ in $T_v(\mathbb{P}T^*N)$ as $$(Tp)_v(\xi_v) = V,$$ where $V$ is the contact element (as hyperplane) associated to $v$. One can check that $\xi: \mathbb{P}T^*N \rightarrow T(\mathbb{P}T^*N): v \mapsto (v,\xi_v)$ defines a contact structure on $\mathbb{P}T^*N$. Notice that the fibers of $p$ (i.e. the contact elements with fixed foot point) form Legendrian submanifolds of $\mathbb{P}T^*N$ (see also exercise \ref{Oefening terugverwijzen 1} of this chapter).\\
Now, let $(M,\xi)$ be a contact manifold, $\mathcal{F}$ a co-Legendrian foliation on $(M,\xi)$ and $\mathcal{G}$ the Legendrian foliation associated to $\mathcal{F}$. Using the theorem of Darboux, we have that $M$ is locally contactomorphic to $\mathbb{P}T^*N$ for an arbitrary smooth manifold $N$ (which we will fix from now on), so we will assume wlog that $M = \mathbb{P}T^*N$ (notice that $\dim(M) = 2\dim(N)-1$). For the same reason, we may assume wlog that the leaves of $\mathcal{G}$ are given by the fibers of $p: \mathbb{P}T^*N \rightarrow N$.\\
The projection of $\mathcal{F}$ under $p$ will give rise to a $1$-dimensional foliation in $N$ which we will call $\mathcal{L}$. Hence a leaf of $\mathcal{F}$ is nothing more than the collection of contact elements for which the foot points all lie in a leaf of $\mathcal{L}$. The leaves of $\mathcal{L}$ are the trajectories of a $1$-parameter group of diffeomorphisms $\Psi^N$ of $N$, explicitely given by $$\Psi^N = \{\Psi_t^N | t \in \mathbb{R}\}.$$ Notice that if $\psi^N$ is a diffeomorphism on $N$, then we get a natural diffeomorphism $\psi^{M}$ on the space of contact elements of $N$, given by $$\psi^{M}(n,V) = (\psi^N(n), (T\psi^N)_n(V)),$$
for $n \in N$ and $V$ a contact element of $N$ in $T_nN$. Using this ``lift'', the $1$-parameter group of diffeomorphisms $\Psi^N$ of $N$ induces a $1$-parameter group of contactomorphisms $\Psi^M$ of $M$ such that the foliation $\mathcal{G}$ is preserved (see also exercise \ref{Oefening terugverwijzen 2} of this chapter). Finally, if we restrict to the open set of contact elements of $N$ which are not tangent to $\mathcal{L}$, the $1$-parameter group $\Psi^M$ allows us to construct $\mathcal{F}$ from $\mathcal{G}$, as we wanted to show.
\end{proof}\noindent
Another more technical but nevertheless crucial property is described in section 2.3 of \cite{Khesin Tabachnikov}. It holds that the leaves of a Legendrian foliation can be endowed with an affine structure (i.e. the leaves of this foliation are locally diffeomorphic to an affine space).
\newline \hfill
\newline
\noindent
\fbox{ \parbox{\textwidth}{
\begin{proposition}
Let $(M,\xi)$ be a $(2n+1)$-dimensional contact manifold with a co-Legendrian foliation $\mathcal{F}$ and associated Legendrian foliation $\mathcal{G}$. The leaves of $\mathcal{G}$ can be endowed with an affine structure.
\end{proposition}}}
\begin{proof}
A sketch of the proof goes as follows: the space of all contact elements in a fixed tangent space is diffeomorphic to the projective space $\mathbb{RP}^{n}$. Using the terminology of lemma \ref{model co legendrian}, we thus have that a leaf $p^{-1}(x)$ of $\mathcal{G}$ can be endowed with a projective structure by the local diffeomorphism $$\varphi: p^{-1}(x) \rightarrow \mathbb{RP}^{n}: y \mapsto (Tp)_y(\xi_v).$$
In Lemma 2.4 in $\cite{Khesin Tabachnikov}$, it has been proven that the image of $\varphi$ actually lies in $\mathbb{RP}^{n}\setminus  \mathbb{RP}^{n-1}$, which is the affine space $\mathbb{A}^{n}$. Hence $\varphi$ defines an affine structure on the leaves of $\mathcal{G}$. For more details, see \cite{Khesin Tabachnikov}.
\end{proof}\noindent
The existence of this affine structure on the leaves of a Legendrian foliation has strong topological and dynamical consequences. On the topological side, the compact leaves of a Legendrian foliation are $n$-dimensional tori.
\newline 
The dynamical consequences of the existence of this affine structure, which will be discussed in the next section, will actually be the fundamental basis on which the (first) notion of contact complete integrability is built.
\section{Contact Completely integrable systems}
As announced at the beginning of this chapter, we will now consider $2$ possible notions for contact complete integrability: the first one is more geometric in nature (and is due to \cite{Khesin Tabachnikov}) and the second one is comparable to the definition that we have seen in the symplectic case (mainly due to \cite{Jovanovic} and \cite{Jovanovic Jovanovic}). Both of them have advantages and disadvantages in their intuition and application.\newline

\noindent
Let $(W,\omega,h)$ be a completely integrable system in the symplectic case. We have seen in chapter $2$ of the lecture notes on symplectic geometry (in particular section $2.1$) that the Liouville foliation of the system (of which the leaves are the connected components of the fibers of $h$) has isotropic leaves. In particular, the regular leaves (consisting of regular points) are Lagrangian submanifolds of $(W,\omega)$. One could have alternatively defined a completely integrable system by considering such a singular Lagrangian foliation (of which almost all leaves are regular) on a symplectic manifold and a vector field $X$ that leaves the foliation invariant (thus the ODE is $z' = X(z)$). The integrals are then (locally) defined as the functions of which the intersection of the level sets form the leaves of the foliation (which you can do by the implicit function theorem).\newline

\noindent
On a $(2n+1)$-dimensional contact manifold $(M,\xi)$, we will also define a contact completely integrable system in such a way. Namely, we consider a contact vector field $X$ that leaves a (singular\footnote{In what follows, we will focus on the regular points, which should form an open and dense subspace in the manifold.}) co-Legendrian foliation $\mathcal{F}$ leaf-wise invariant. Since $X$ is a contact vector field, the leaves of the associated Legendrian foliation $\mathcal{G}$ within a leaf $F$ of $\mathcal{F}$ will be send to each other under the flow of $X$. Furthermore, the affine structure on the leaves of $\mathcal{G}$ will be preserved by the flow. Hence the only thing we still need to know about the dynamics, is how the leaves of $\mathcal{G}$ are send to each other. So we need to consider the dynamics on the $1$-dimensional quotient $F/\mathcal{G}$. \newline

\noindent
A natural notion of integrability on a $1$-dimension manifold $L$, is the existence of an invariant $1$-form $\beta$ on $L$. Indeed, $\beta$ is a closed $1$-form on $L$ (it is top-dimensional). Hence by Poincar\'e's Lemma, $\beta$ can locally be described by an exact $1$-form $df$ for a function $f$ on $L$ known as the integral.\newline 

\noindent
Using all of these insights, we will first define a completely integrable contact manifold (a purely geometric construction which includes the foliation) and then define a contact completely integrable system as a contact vector field which satisfies some properties.
\newline \hfill
\newline
\noindent
\fbox{ \parbox{\textwidth}{
\begin{definition}\label{def 1}
The triple $(M,\xi, \mathcal{F})$ is called a $(2n+1)$-dimensional \textit{completely integrable contact manifold} if:
\begin{itemize}
\item $(M,\xi)$ is a $(2n+1)$-dimensional contact manifold,
\item $\mathcal{F}$ is a co-Legendrian foliation on $(M,\xi)$,
\item for each leaf $F$ of $\mathcal{F}$, the associated Legendrian foliation $\mathcal{G}|_F$ has a holonomy invariant transverse smooth $1$-form.
\end{itemize}
Let $(M,\xi,\mathcal{F})$ be a completely integrable contact manifold. A vector field $X$ on $M$ is called a \textit{contact completely integrable system} on $(M,\xi,\mathcal{F})$ if
\begin{itemize}
\item $X$ is a contact vector field on $(M,\xi)$.
\item $X$ is tangent to the co-Legendrian foliation $\mathcal{F}$ (so it leaves $\mathcal{F}$ leaf-wise invariant).
\item $X$ preserves the holonomy invariant transverse smooth $1$-forms.
\end{itemize}
\end{definition}}}\newline \hfill \newline \noindent
The notion of a holonomy invariant transverse smooth $1$-form is something technical from advanced differential geometry\footnote{Holonomy will in some way be a generalization of monodromy, which we will encounter at the end of the lecture notes on symplectic geometry in the context of semitoric systems.}. In essence, these conditions ensure that the $1$-form will descend to a $1$-form on the $1$-dimensional quotient $F/\mathcal{G}$. For more details, we refer to \cite{Khesin Tabachnikov} and the references therein.\\
Luckily, it is possible to completely ignore these technical details if the reader prefers to. Indeed, a (slight adaptation of a) characterisation described in \cite{Khesin Tabachnikov} reformulates the definition above such that the technical stuff is hidden in the background.
\newline \hfill
\newline
\noindent
\fbox{ \parbox{\textwidth}{
\begin{theorem}\label{theorem def 1}
The triple $(M,\xi, \mathcal{F})$ is a $(2n+1)$-dimensional completely integrable contact manifold if and only if:
\begin{itemize}
\item $(M,\xi)$ is a $(2n+1)$-dimensional contact manifold,
\item $\mathcal{F}$ is a co-Legendrian foliation on $(M,\xi)$,
\item there exists a contact form $\alpha$ on $M$ such that (locally) $\xi = \ker(\alpha)$ and the Reeb vector field $R_\alpha$ is tangent to $\mathcal{F}$.
\end{itemize}
The contact form $\alpha$ as described above is called the \textit{contact form on $M$ that defines the completely integrable contact manifold} $(M,\xi,\mathcal{F})$.\newline

\noindent
Let $(M,\xi,\mathcal{F})$ be a completely integrable contact manifold with defining contact form $\alpha$. A vector field $X$ on $M$ is a contact completely integrable system on $(M,\xi,\mathcal{F})$ if and only if
\begin{itemize}
\item $X$ is tangent to the co-Legendrian foliation $\mathcal{F}$.
\item $X$ is a strict contact vector field with respect to $\alpha$.
\end{itemize}
\end{theorem}}}
\begin{proof}
This result is one of the main topics in \cite{Khesin Tabachnikov}.
\end{proof}\noindent
If $(M,\xi,\mathcal{F})$ is a completely integrable contact manifold with defining contact form $\alpha$, then the restriction of $\alpha$ to a leaf of $\mathcal{F}$ will define the holonomy invariant transverse smooth measure. This characterisation should be more easy to understand (you don't need advanced differential geometry). However, the downside is that it is less clear what the intuition is behind the definition in comparison to the original definition (which is the reason why we have chosen to introduce it in this order).\newline

\noindent
There exists a second notion of complete integrability in the contact setting (there exist even more notions which all have their advantages and disadvantages, but we will restrict ourselves to $2$), as introduced by \cite{Jovanovic} and \cite{Jovanovic Jovanovic}\footnote{They even introduced non-commutative contact complete integrability, of which our definition is a special case.}. This definition will only work on strict contact manifolds $(M,\alpha)$.
\newline \hfill
\newline
\noindent
\fbox{ \parbox{\textwidth}{
\begin{definition}\label{def 2}
Let $(M,\alpha)$ be a $(2n+1)$-dimensional strict contact manifold and consider the differential equation
\begin{equation*}
z' = X(z)
\end{equation*}
for a contact vector field $X$. The equation is called contact completely integrable if there exist smooth functions $f_1,\cdots, f_n$ on $M$ such that
\begin{itemize}
\item $f_i$ is an integral of $X$ for $i = 1,\cdots, n$ (i.e. $X(f_i) = 0$).
\item The functions will pairwise commute with respect to the Jacobi bracket $$[f_i,f_j] = 0 \qquad \textrm{for } i= 1,\cdots, n \textrm{ and } j = 1,\cdots,n.$$
\item The $1$-forms $df_1,\cdots,df_n$ are almost everywhere linearly independent.
\item $\alpha(X) = f_1$ (so $X$ is the contact Hamiltonian of $f_1$) or $\alpha(X) = 1$ (so $X = R_\alpha$).
\item The functions  $f_1,\cdots,f_n$ will commute with the constant function $\mathbbm{1}$ with respect to the Jacobi bracket $$[\mathbbm{1},f_i] = 0 \qquad \textrm{for } i= 1,\cdots, n.$$
\end{itemize}
\end{definition}}}
In the rest of this chapter, we will show that both notions of contact complete integrability are (more or less) equivalent. We will work on a strict contact manifold.\newline

\noindent
Consider the (global) contact form $\alpha$ on a $(2n+1)$-dimensional manifold $M$ defining a completely integrable contact manifold $(M,\xi = \ker \alpha, \mathcal{F}$). Let $X$ be a contact completely integrable system on $(M,\xi,\mathcal{F})$.  Using the implicit function theorem, there exist (locally but we can patch them together at the end) $n$ functions $f_1,\cdots,f_n$ such that the leaves of $\mathcal{F}$ are given by the intersection of the level sets of $f_1,\cdots,f_n$. Since $X$ is tangent to $\mathcal{F}$, it holds that $f_1,\cdots,f_n$ are integrals of $X$ (can you see why?). We will now prove a series of lemmas which we will need in order to show that the system is contact completely integrable in the sense of definition \ref{def 2}.
\newline \hfill
\newline
\noindent
\fbox{ \parbox{\textwidth}{
\begin{lemma}
Let $f$ be an integral of $\mathcal{F}$ (so $f$ is constant on the leaves of $\mathcal{F}$). It holds that $X_f$ is tangent to $\mathcal{F}$.
\end{lemma}}}
\begin{proof}
Let $f$ be an integral of $\mathcal{F}$ and let $\mathcal{G}$ be the associated Legendrian foliation to $\mathcal{F}$. Let $w \in T\mathcal{G}$ be arbitrary. It holds that
\begin{equation*}
d\alpha(X_f,w) = df(R_\alpha) \alpha(w) - df(w).
\end{equation*}
Since $T\mathcal{G} \subset \xi$, we have that $\alpha(w) = 0$. Since $f$ is an integral of $\mathcal{F}$ and since $w \in T\mathcal{F}$, it holds that $df(w) = 0$. Hence $d\alpha(\widetilde{X_f},w) = d\alpha(X_f,w) = 0$, so $\widetilde{X_f}$ lies in the symplectic complement of $T\mathcal{G}$ in $\xi$. However, since $\mathcal{G}$ is Legendrian, this means (by a dimension count) that $\widetilde{X_f} \in T\mathcal{G}$. Since it also holds that $R_\alpha \in T\mathcal{F}$, we conclude that
\begin{equation*}
X_f = fR_\alpha + \widetilde{X_f} \in T\mathcal{F}.
\end{equation*}
This means exactly that $X_f$ is tangent to $\mathcal{F}$.
\end{proof}
\noindent
\fbox{ \parbox{\textwidth}{
\begin{lemma}
Consider the situation described above. The functions $f_1,\cdots,f_n$ will pairwise commute with respect to the Jacobi bracket.
\end{lemma}}}
\begin{proof}
Let $i,j \in \{1,\cdots,n\}$. It holds that $df_i(R_\alpha) = df_j(R_\alpha) = 0$ since $R_\alpha \in T\mathcal{F}$. Hence we find that
\begin{eqnarray*}
[f_i,f_j] &=& d\alpha(X_{f_i},X_{f_j}) + f_i R_\alpha(f_j) - f_j R_\alpha(f_i)\\
&=& d\alpha(X_{f_i},X_{f_j}) \\
&=& d\alpha(\widetilde{X_{f_i}},\widetilde{X_{f_j}}) \\
&=& 0,
\end{eqnarray*}
where the last equality holds since $\widetilde{X_{f_i}},\widetilde{X_{f_j}} \in T\mathcal{G}$ (we use the previous lemma here) and $\mathcal{G}$ is isotropic with respect to $(\xi,d\alpha|_\xi)$.
\end{proof}
\noindent
\fbox{ \parbox{\textwidth}{
\begin{lemma}
Consider the situation described above. The functions $f_1,\cdots,f_n$ will commute with the constant function $\mathbbm{1}$ with respect to the Jacobi bracket.
\end{lemma}}}
\begin{proof}
Let $i \in \{1,\cdots,n\}$.  We have again that $df_i(R_\alpha) = 0$. Consequently, it holds that
\begin{eqnarray*}
[\mathbbm{1},f_i] &=& d\alpha(R_\alpha,X_{f_i}) + \mathbbm{1} R_\alpha(f_i) - f_j R_\alpha(\mathbbm{1})\\
&=& d\alpha(R_\alpha,X_{f_i}) \\
&=& 0,
\end{eqnarray*}
so we are done.
\end{proof}
\noindent
\fbox{ \parbox{\textwidth}{
\begin{lemma}
Consider the situation as it is described above. It holds that $df_1,\cdots,df_n$ are (almost everywhere) linear independent.
\end{lemma}}}
\begin{proof}
$\mathcal{F}$ is in fact a singular co-Legendrian foliation such that almost all points lie on a leaf of dimension $n+1$ of $\mathcal{F}$ (but in the earlier discussion, we have ignored these singular points and assumed for the sake of notation that every point is regular). Consider a leaf $F$ of dimension $n+1$ and let $x \in F$. Since $f_1,\cdots,f_n$ determine $F$, we have that $$(T_xF)^\circ = \textrm{span}\{(df_1)_x,\cdots,(df_n)_x\}.$$ Furthermore, since $\dim((T_xF)^\circ) = \dim(M)-\dim(F) = n$, it holds that $(df_1)_x,\cdots,(df_n)_x$ are linear independent . Hence $df_1,\cdots,df_n$ are almost everywhere linear independent.
\end{proof}
\noindent
\fbox{ \parbox{\textwidth}{
\begin{lemma}\label{lemma lin onaf}
Consider the situation as it is described above. It holds that $R_\alpha,X_{f_1},\cdots,X_{f_n}$ are (almost everywhere) linear independent.
\end{lemma}}}
\begin{proof}
Consider a leaf of $\mathcal{F}$ of dimension $n+1$. By the previous lemma, $df_1,\cdots,df_n$ are linear independent over this leaf. Now, let $c_0,\cdots,c_n \in \mathbb{R}$ such that
\begin{equation}\label{verg lin onaf}
c_0 R_\alpha + \sum_{i=1}^n c_iX_{f_i} = 0
\end{equation}
If we evaluate this vector in $d\alpha$, we get the $1$-form
\begin{eqnarray*}
0 &=& d\alpha(c_0 R_\alpha + \sum_{i=1}^n c_iX_{f_i},\cdot) \\
&=& \sum_{i=1}^n c_i d\alpha(X_{f_i},\cdot)\\
&=&  \sum_{i=1}^n c_i (df_i(R_\alpha)\alpha - df_i)\\
&=&  -\sum_{i=1}^n c_i df_i,
\end{eqnarray*}
where we have used in the last equation that $R_\alpha$ is tangent to $\mathcal{F}$. Using the linear independence of $df_1,\cdots,df_n$ over the leaf, we find that $c_i = 0$ for $ i = 1,\cdots,n$. Hence equation (\ref{verg lin onaf}) becomes $c_0 R_\alpha = 0$, so $c_0 = 0$. We conclude that $R_\alpha,X_{f_1},\cdots,X_{f_n}$ are (almost everywhere) linear independent.
\end{proof}
\noindent
\fbox{ \parbox{\textwidth}{
\begin{corollary}
Consider the situation as described above. The vector fields $R_\alpha,X_{f_1},\cdots,X_{f_n}$ will span the foliation $\mathcal{F}$.
\end{corollary}}}
\begin{proof}
If $F$ is a regular leaf of $\mathcal{F}$ (so $\dim(F) = n+1$), we have that $R_\alpha,X_{f_1},\cdots,X_{f_n}$ are linear independent in $TF$ and the dimension of a tangent space of $F$ is $n+1$, so we are done. The same argument holds for singular leaves, but you will find less integrals that describe the leaf, so less vector fields that span the tangent bundle.
\end{proof}
\noindent
\fbox{ \parbox{\textwidth}{
\begin{lemma}
Consider the situation as described above. If $X$ is not (a rescaling of) $R_\alpha$, then we can say wlog that $f_1 = \alpha(X)$.
\end{lemma}}}
\begin{proof}
It holds that $X$ is a rescaling of $R_\alpha$ if and only if $\alpha(X)$ is a constant function.  Hence $\alpha(X)$ is not a constant function. Assume that we have integrals $f_1,\cdots,f_n$ as described above. Using Cartan's magic formula and the fact that $X$ is a strict contact vector field with respect to $\alpha$ (equivalently $\mathcal{L}_X\alpha = 0$), we have that
\begin{eqnarray*}
R_\alpha(\alpha(X)) &=& d(\iota_X \alpha) (R_\alpha)\\
&=& -(\iota_Xd\alpha) (R_\alpha) + (\mathcal{L}_X\alpha)(R_\alpha)\\
&=& -d\alpha (X,R_\alpha)\\
&=& 0.
\end{eqnarray*}
If $i \in \{1,\cdots,n\}$, we will find for the same reason (and since $df_i(R_\alpha)$ and $df_i(X)$ are both $0$) that
\begin{eqnarray*}
X_{f_i}(\alpha(X)) &=& d(\iota_X \alpha) (X_{f_i})\\
&=& -(\iota_Xd\alpha) (X_{f_i}) + (\mathcal{L}_X\alpha)(X_{f_i})\\
&=& -d\alpha (X,X_{f_i})\\
&=& d\alpha (X_{f_i},X)\\
&=& df_i(R_\alpha)\alpha(X) - df_i(X)\\
&=& 0.
\end{eqnarray*}
This means that $\alpha(X)$ is an integral of the foliation $\mathcal{F}$. Since it is not constant as a function, we know that $d(\alpha(X)) \neq 0$, so $d(\alpha(X))$ is a non-zero element in the annihilator of the leaves. This means that we can construct a basis of the annihilator of the leaves which includes $\alpha(X)$. This basis gives you the other integrals, so wlog $f_1 = \alpha(X)$.
\end{proof}\noindent
If we put all these lemmas together, we get exactly that the ODE $z' = X(z)$ is contact completely integrable (with integrals $f_1,\cdots,f_n$) as described in definition \ref{def 2}.\newline

\noindent
Conversely, consider now a $(2n+1)$-dimensional strict contact manifold $(M,\alpha)$ (and denote $\xi = \ker \alpha$). Let $X$ be a contact vector field such that the associated ODE is contact completely integrable with respect to the functions $f_1,\cdots,f_n$ on $M$. Define the leaves of a foliation $\mathcal{F}$ as the intersections of the level sets of $f_1,\cdots,f_n$.
\newline \hfill
\newline
\noindent
\fbox{ \parbox{\textwidth}{
\begin{lemma}
Consider the situation described above. The Reeb vector field $R_\alpha$ associated to $(M,\alpha)$ is tangent to the foliation $\mathcal{F}$.
\end{lemma}}}
\begin{proof}
Let $F$ be a leaf of $\mathcal{F}$. By construction of the foliation $\mathcal{F}$, we know that $F$ is the intersection of certain level sets of the functions $f_1,\cdots,f_n$. Let $j \in \left\{1,\cdots, n \right\}$. Since the function $f_j$ satisfies the conditions as described in definition \ref{def 2}, we know that $\left[\mathbbm{1},f_j\right] = 0$, which means that $R_\alpha(f_j) = 0$. If $\varphi_t$ is the local flow map of the  vector field $R_\alpha$, we have for every $z \in F$ that
$$0 = R_\alpha(f_j)(z) = df_j|_z(R_\alpha(z)) = \frac{d}{dt}\Big|_{t=0} f_j\left(\varphi_t(z)\right).$$
Consequently, the flow of $R_\alpha$ stays in a level set of $f_j$. Since this holds for every $j$, the flow of $R_\alpha$ will stay in the intersections of the level sets of the functions, so it will stay in $F$. Hence $R_\alpha$ is tangent to $F$.
\end{proof}
\noindent
\fbox{ \parbox{\textwidth}{
\begin{lemma}
Consider the situation described above. The (singular) foliation $\mathcal{F}$ is a (singular) co-Legendrian foliation on $(M,\xi)$.
\end{lemma}}}
\begin{proof}
Since $df_1,\cdots,df_n$ are almost everywhere linear independent, we have that a generic leaf of $\mathcal{F}$ is $(n+1)$-dimensional. From now on, we will work over such a generic leaf $F$. The contact Hamiltonian $X_{f_i}$ will be tangent to the level set of $f_j$. Indeed, if $\varphi_i$ is the flow of $X_{f_i}$, it holds that $X_{f_i}$ is tangent to the level set of $f_j$ if and only if $(f_j \circ \varphi_i)' = 0$. However
\begin{equation*}
(f_j \circ \varphi_i)' = df_j(X_{f_i}) = X_{f_i}(f_j).
\end{equation*}
Since $[1,f_i] = [1,f_j] = 0$, we have by Proposition \ref{Jacobi bracket als Poisson bracket} that $X_{f_i}(f_j) = 0$ if and only if $[f_i,f_j] = 0$. We conclude that $X_{f_i}$ is tangent to $F$. Furthermore, since $R_\alpha$ is tangent to $F$, we have that $$\textrm{span}(R_\alpha,X_{f_1},\cdots,X_{f_n}) \subseteq TF.$$
By copying the proof of lemma \ref{lemma lin onaf}, we have that $R_\alpha,X_{f_1},\cdots,X_{f_n}$ are linear independent, so both sides of the inclusion above are $(n+1)$-dimensional vector spaces (when evaluated in a point), which gives us that $$\textrm{span}(R_\alpha,X_{f_1},\cdots,X_{f_n}) = TF.$$
Now, consider the distribution $TF\cap \xi$. Since both $R_\alpha$ and $X_{f_i}$ are in $TF$, it holds that the horizontal part $\widetilde{X_{f_i}}$ is in $TF$, so $\widetilde{X_{f_i}} \in TF\cap \xi$. Consequently $$\textrm{span}(\widetilde{X_{f_1}},\cdots,\widetilde{X_{f_n}}) \subseteq TF\cap \xi.$$
On the one hand, $\widetilde{X_{f_1}},\cdots,\widetilde{X_{f_n}}$ are linearly independent (since so are $R_\alpha,X_{f_1},\cdots,X_{f_n}$). On the other hand, $\mathcal{F}$ is transverse to $\xi$ (since $R_\alpha \in T\mathcal{F}$), so $\dim(TF\cap \xi) = n$. We conclude that $$\textrm{span}(\widetilde{X_{f_1}},\cdots,\widetilde{X_{f_n}}) = TF\cap \xi.$$ It follows from $[f_i,f_j] = 0$ and $[1,f_i] = 0$ that $[X_{f_i},X_{f_j}]=0$ and $[R_\alpha,X_{f_i}] = 0$. Hence by a standard computation of Lie brackets of vector fields, we find that
\begin{equation*}
[\widetilde{X_{f_i}},\widetilde{X_{f_j}}] = [X_{f_i} - f_iR_\alpha,X_{f_j} - f_jR_\alpha] = 0.
\end{equation*}
This means that $TF\cap\xi$ is an involutive distribution. Hence by the Theorem of Frobenius, it is an integrable distribution. Putting everything together, we have that $\mathcal{F}$ is a co-Legendrian foliation.
\end{proof}
\noindent
\fbox{ \parbox{\textwidth}{
\begin{corollary}
Consider the situation described above. It holds (almost everywhere) that $$\textrm{span}(R_\alpha,X_{f_1},\cdots,X_{f_n}) = T\mathcal{F}.$$
\end{corollary}}}
\begin{proof}
This result has been proven in the proof of the previous lemma.
\end{proof}
\noindent
\fbox{ \parbox{\textwidth}{
\begin{corollary}
Consider the situation described above. The vector field $X$ is tangent to $\mathcal{F}$.
\end{corollary}}}
\begin{proof}
If we consider the fourth condition in definition \ref{def 2}, we find that $X \in \{X_{f_1},R_\alpha\}$. Consequently, we have that $$X \in \textrm{span}(R_\alpha,X_{f_1},\cdots,X_{f_n}),$$
so $X$ is tangent to $\mathcal{F}$.
\end{proof}
\noindent
\fbox{ \parbox{\textwidth}{
\begin{lemma}
Consider the situation described above. The vector field $X$ is a strict contact vector field with respect to $\alpha$.
\end{lemma}}}
\begin{proof}
Again, we have that $X \in \{X_{f_1},R_\alpha\}$. If $X = R_\alpha$, then we are done by example \ref{Reeb VF strict contact VF}. Consider the case where $X = X_{f_1}$. We find, using Cartan's magic formula, that
\begin{equation*}
\mathcal{L}_{X_{f_1}}\alpha = d(\alpha(X_{f_1})) + \iota_{X_{f_1}}d\alpha = df_1 + df_1(R_\alpha)\alpha - df_1 = df_1(R_\alpha)\alpha.
\end{equation*}
Since we know that $[1,f_1] = 0$, it holds that $\mathcal{L}_{X_{f_1}}\alpha = df_1(R_\alpha)\alpha = 0$. We conclude that $X_{f_1}$ is a strict contact vector field with respect to $\alpha$ as well.
\end{proof}\noindent
If we now put everything together, we have proven that $X$ defines a contact completely integrable system on the completely integrable contact manifold $(M,\xi,\mathcal{F})$.\newline

\noindent
We can summarize the discussions above as follows.
\newline \hfill
\newline
\noindent
\fbox{ \parbox{\textwidth}{
\begin{theorem}
Contact complete integrability in the sense of definition \ref{def 1} and theorem \ref{theorem def 1} is equivalent to contact complete integrability of the Reeb type in the sense of definition \ref{def 2} and definition \ref{def Reeb type}.
\end{theorem}}}\newline \hfill \newline\noindent
We will refer to the exercises of this section for a concrete (non-trivial) example of a contact completely integrable system. 
\newpage
\section{Exercises}

\begin{enumerate} 
\item Consider $\mathbb{R}^3$ endowed with Cartesian coordinates $(x,y,z)$. Furthermore, consider $\alpha_{st} = dz + xdy$, the standard contact form on $\mathbb{R}^3$.
\begin{enumerate}
\item Find a maximal open subset $U \subset \mathbb{R}^3$ such that the distribution $$\mathcal{F}_1 = \textrm{span}\left(\frac{\partial}{\partial x},\frac{\partial}{\partial x}-y\frac{\partial}{\partial z}\right)$$ defines a co-Legendrian foliation on $(U,\alpha_{st}|_U)$. 
\item Describe the Legendrian foliation associated to this co-Legendrian foliation.
\item Decide whether $(U,\ker \alpha_{st}|_U,\mathcal{F}_1)$ is a completely integrable contact manifold.
\end{enumerate}

\item Consider $\mathbb{R}^3$ endowed with Cartesian coordinates $(x,y,z)$. Furthermore, consider $\alpha_{st} = dz + xdy$, the standard contact form on $\mathbb{R}^3$.
\begin{enumerate}
\item Show that the foliation defined by the (integrable) distribution $$\mathcal{F}_2 = \textrm{span}\left(\frac{\partial}{\partial x},\frac{\partial}{\partial z}\right)$$ is a co-Legendrian foliation on $(\mathbb{R}^3,\alpha_{st})$. 
\item Describe the Legendrian foliation associated to this co-Legendrian foliation.
\item Show that $(\mathbb{R}^3,\ker \alpha_{st},\mathcal{F}_2)$ is a completely integrable contact manifold.
\item Find a smooth function $f$ on $\mathbb{R}^3$ such that the foliation $\mathcal{F}_2$ can be described as the collection of level sets of $f$.
\item Show that the contact Hamiltonian vector field $X_f$ is a contact completely integrable system on $(\mathbb{R}^3,\ker \alpha_{st},\mathcal{F}_2)$.
\end{enumerate}

\item  In this exercise, we will take a closer look at the proof of proposition \ref{model co legendrian}. We will refer to this proof for the notations and context of this exercise.
\begin{enumerate}
\item \label{Oefening terugverwijzen 1} Show that the fibers of $p$ are Legendrian submanifolds of $\mathbb{P}T^*N$.
\item \label{Oefening terugverwijzen 2} Show that the lift of the $1$-parameter group of diffeomorphisms $\Psi^N$ of $N$ induces a $1$-parameter group of contactomorphisms $\Psi^M$ of $M$ such that the foliation $\mathcal{G}$ is preserved (which means that leaves of $\mathcal{G}$ are send to leaves of $\mathcal{G}$).
\end{enumerate}

\item Consider the manifold $\mathbb{C}^{n+1}$, which can be identified with $\mathbb{R}^{2n+2}$ via the identification $z_j = x_j + i y_j$ for $j = 0,\cdots,n$. Furthermore, consider the smooth function
$$F(z,\overline{z}) = \sum_{j=0}^n |z_j|^2$$
on $\mathbb{C}^{n+1}$. Then the unit sphere $\mathbb{S}^{2n+1} = F^{-1}(1)$ is a subset of $\mathbb{C}^{n+1}$.\\
Now, let $a_j$ be a positive real number for $j = 0,\cdots,n$ and consider the $1$-form
$$\alpha = \frac{i}{8}\sum_{j=0}^n a_j(z_jd\overline{z}_j - \overline{z}_jdz_j).$$
One can show (you don't have to do this) that $(\mathbb{S}^{2n+1},\alpha)$ is a strict contact manifold with Reeb vector field
$$R_\alpha = 4i \sum_{j=0}^n \frac{1}{a_j}\left(z_j\frac{\partial}{\partial z_j} - \overline{z}_j \frac{\partial}{\partial \overline{z}_j} \right).$$
Finally, we will also need the functions
$$\forall j \in \{0,\cdots,n\}: f_j(z) = |z_j|^2$$
and the vector fields on $\mathbb{C}^{n+1}$ given by
$$\forall j \in \{0,\cdots,n\}: Y_j = \frac{4i}{a_j}\left(z_j\frac{\partial}{\partial z_j} - \overline{z}_j \frac{\partial}{\partial \overline{z}_j} \right).$$
\begin{enumerate}
\item Show that $f_j$ is an integral of $R_\alpha$ for every $j \in \{0,\cdots,n\}$.\\
\item Show that $Y_j$ is tangent to $\mathbb{S}^{2n+1}$ for every $j \in \{0,\cdots,n\}$.
\item Show that $Y_j$ is the contact Hamiltonian vector field of $f_j$ for every $j \in \{0,\cdots,n\}$.
\item Show that the functions $f_j$ and $f_k$ commute with respect to the Jacobi bracket on $(\mathbb{S}^{2n+1},\alpha)$ for every $j,k \in \{0,\cdots,n\}$.
\item Show that there exists an open dense subset $U$ of $\mathbb{C}^{n+1}$ such that it holds on $U$ that
$$dF \wedge df_1 \wedge \cdots \wedge df_n|_U \neq 0|_U.$$
One can conclude from this (you don't have to do this) that 
$$df_1 \wedge \cdots \wedge df_n|_{\mathbb{S}^{2n+1}\cap U} \neq 0|_{\mathbb{S}^{2n+1}\cap U}.$$
\textit{\underline{Hint:} Work with polar coordinates $z_j = r_je^{i\varphi_j}$.}
\item Show that the equation $$z' = R_\alpha(z)$$ defines a contact completely integrable system on $(\mathbb{S}^{2n+1},\alpha)$ with integrals $f_1,\cdots,f_n$ in the sense of Definition \ref{def 2}.
\end{enumerate}

\item Explain why the existence of an affine structure on the leaves of a Legendrian foliation implies that the compact leaves of this foliation are tori.

\item Compare this definition of contact complete integrability with the definition of complete integrability in the symplectic case. Discuss in particular the similarities and the differences.

\end{enumerate}

\chapter{(Semi-)Local aspects}
In this (relatively short) chapter, we will study some local (in fact even semi-local) aspects of contact complete integrability. In particular, we will give a contact variant of the Arnold-Liouville Theorem (by \cite{Jovanovic}) and discuss a local normal form in the neighbourhood of certain critical points (by \cite{Miranda}). Morally, this chapter coincides with the second chapter of the lecture notes on symplectic geometry. The main goal is to give the reader a taste of these results, but we refer to the given sources for the technical details and proofs (which are comparable in size to the proofs of chapter $2$ of the lecture notes symplectic geometry). In this chapter, we will use the concept of contact complete integrability given by definition \ref{def 2}.\newline

\noindent
Just like in the symplectic case, the contact version of the Arnold-Louville Theorem has a topological part (in which we will describe the topology of the compact connected components of certain level sets) and a dynamical part (in which we will introduce some sort of action-angle variables and in which we will see that the flow on the level sets above is quasi-periodic).
\newline \hfill
\newline
\noindent
\fbox{ \parbox{\textwidth}{
\begin{theorem}\label{Contact AL}
\textbf{(Arnold-Liouville)} Let $(M,\alpha)$ be a $(2n+1)$-dimensional strict contact manifold and assume that the differential equation $$z' = X_f(z)$$ is contact completely integrable (in the sense of definition \ref{def 2}) with integrals $f_1,\cdots,f_n$ (in particular, it holds that $f=1$ or $f=f_1$). Let $F$ be a compact connected component of the level set $$\{x \in M | f_i(x) = c_i \quad \textrm{ for } i = 1,\cdots,n\}$$ on which it holds that $$df_1\wedge\cdots\wedge df_n \neq 0.$$
It holds that
\begin{itemize}
\item[a)] $F$ is diffeomorphic to the $(n+1)$-dimensional torus $\mathbb{T}^{n+1}$.
\item[b)] There exists an open neighbourhood $U \subset M$ of $F$ and a contact form $\alpha$ on $U$, both invariant under the Hamiltonian flows of the integrals $f_j$, and there exists a diffeomorphism $\varphi: U \rightarrow \mathbb{T}^{n+1}\times D$, given by $$\forall m \in U: \varphi(m) = (\theta,y) = (\theta_0,\cdots,\theta_n,y_1,\cdots,y_n),$$
for a certain open set $D \subset \mathbb{R}^n$ such that
\begin{itemize}
\item[i)] The contact form $\alpha$ takes on the form of the following contact form $\alpha_0$ on $\mathbb{T}^{n+1}\times D$:
\begin{equation*}
\alpha_0 = (\varphi^{-1})^*\alpha = y_0d\theta_0 + y_1d\theta_1 + \cdots + y_nd\theta_n,
\end{equation*}
where $y_0$ is a smooth function on $\mathbb{T}^{n+1}\times D$ that only depends on $y_1,\cdots,y_n$.
\item[ii)] The flow of the contact Hamiltonian $X_f$ on the invariant tori (such as $F$) will be quasi-periodic. This means that it is of the form
\begin{equation*}
\forall t \in \mathbb{R}: (\theta_0,\cdots,\theta_n) \mapsto (\theta_0 + t\omega_0,\cdots,\theta_n + t\omega_n), 
\end{equation*}
where the terms $\omega_0,\cdots,\omega_n$ are called the frequencies, which will only depend on $y_1,\cdots,y_n$.
\end{itemize}
\end{itemize}
\end{theorem}}}
\newpage
\begin{proof}
This theorem is a direct consequence of Theorem $5.1$ and Theorem $5.2$ in \cite{Jovanovic}. In fact, these theorems state an even stronger result, namely that it is only necessary for the differential equation to be noncommutatively contact completely integrable. The proof is based on the construction of a special ``local bi-fibration''.
\end{proof}\noindent
We will call the coordinates $(\theta,y)$ as described in the theorem above the \textit{contact action-angle variables}.
\newline \hfill
\newline
\noindent
\fbox{ \parbox{\textwidth}{
\begin{question} \label{Question 7}
Compare theorem \ref{Contact AL} to the Arnold-Liouville Theorem  in the symplectic case: what are the most noticeable similarities and differences? In particular, compare the role of (symplectic) action-angle variables to contact action-angle variables. Consider also the assumptions in both theorems.
\end{question}}}\newline \hfill \newline \noindent
On the other hand, we can consider the behaviour of a contact completely integrable system near singular points. The following  description can be found in \cite{Miranda} and is, up to the authors knowledge, the ``best" result regarding behaviour near singular points in the contact setting to this date. Even though it is the strongest result so far, we will have to make some additional assumptions which can hopefully be generalized in the near future.\newline

\noindent
Consider a $(2n+1)$-dimensional strict contact manifold $(M,\alpha)$. The first assumption is that we consider the differential equation associated to the Reeb vector field (which was the case where $\alpha(X)=1$) $$z' = R_\alpha(z).$$
Let this equation be contact completely integrable with respect to some integrals $f_1,\cdots,f_n$. Furthermore, we will make the following assumptions:
\begin{itemize}
\item[i)] The Reeb vector field $R_\alpha$ is the infinitesimal generator of an $\mathbb{S}^1$-action, i.e. the Reeb flow will be periodic.
\item[ii)] Consider the matrix $\left(df_1,\cdots,df_n\right)$. Denote the minimal rank\footnote{The rank of a point is still the same as in the symplectic case.} of this matrix as $k$. Let $p\in M$ be a point where the rank is $k$ and denote the orbit through $p$ of the action induced by $X_{f_1},\cdots,X_{f_n}$ as $\mathcal{O}$. We assume that $\mathcal{O}$ is diffeomorphic to a $(k+1)$-dimensional torus $ \mathbb{T}^{k+1}$. If $k<n$, then the point is called singular.
\item[iii)] We assume that the integrals $f_1,\cdots, f_k$ are regular along $\mathcal{O}$. Furthermore the singularities of $f_{k+1},\cdots,f_n$ are non-degenerate in the Morse-Bott sense.\footnote{If the singularity is just a point $q$, then this means that $df_j(q) = 0$ and the Hessian $d^2f_j(q)$ is non-degenerate as a bilinear form for all $j = k+1,\cdots,n$. If the critical points are not isolated but appear as non-zero dimensional submanifolds, then a generalization of this notion of non-degeneracy is possible. We will not go into more details in these notes and refer the interested reader to the literature (search term: Morse-Bott function).}
\end{itemize}
\noindent
\fbox{ \parbox{\textwidth}{
\begin{theorem}\label{Theorem Singular}
\textbf{(Local normal form)} Consider the contact completely integrable system $z' = R_\alpha(z)$ with the additional assumptions and notations described above. In a finite covering of a (tubular) neighbourhood of the orbit $\mathcal{O}$, we can find coordinates $(\theta_0,\cdots,\theta_k,p_1,\cdots,p_k,x_1,y_1,\cdots,x_{n-k},y_{n-k})$ which will take value in $\mathbb{T}^{k+1}\times U^{k}\times V^{2(n-k)}$ (with $U^{k}$ resp. $V^{2(n-k)}$ a $k$-dimensional resp. $2(n-k)$-dimensional disk) such that:
\begin{itemize}
\item[i)] The Reeb vector field in these coordinates is $R_\alpha = \frac{\partial}{\partial\theta_0}$.
\item[ii)] The integrals take on  a special form in these coordinates. For the first $k$ integrals, we have that
\begin{equation*}
f_i = p_i \quad \textrm{for } i = 1,\cdots,k.
\end{equation*}
For the other integrals (so for $j = k+1,\cdots, n$), we get something from the following list:
\begin{itemize}
\item[a)]\textbf{Elliptic component:} $f_j = x_j^2 + y_j^2$,
\item[b)]\textbf{Hyperbolic component:} $f_j = x_jy_j$,
\item[c)]\textbf{Focus-focus component:} These components always appear as a pair ($f_j,f_{j+1}$) for $k < j < n$, namely
\begin{equation*}
\begin{cases}
f_j &= x_jy_{j+1}-x_{j+1}y_j \quad \textrm{and}\\
f_{j+1} &= x_jy_j + x_{j+1}y_{j+1}
\end{cases}
\end{equation*}
\end{itemize}
\item[iii)] In these coordinates, the contact form can be written as $$\alpha_0 = d\theta_0 + \sum_{i=1}^k p_id\theta_i + \frac{1}{2}\sum_{i=1}^{n-k}x_idy_i-y_idx_i.$$
\end{itemize}
The pair $(\mathbb{T}^{k+1}\times U^{k}\times V^{2(n-k)},\alpha_0)$ will be called the contact model manifold.
\end{theorem}}}
\begin{proof}
See Theorem $7.3.1$ in \cite{Miranda}.
\end{proof}
\begin{remark}
Consider the special case where $k=n$. This means, due to the minimality of $k$, that all points of $M$ are regular with respect to this integrable system. In this case, we see that the results of theorem \ref{Contact AL} and theorem \ref{Theorem Singular} are very similar (consider $(\theta,p)$ as contact action-angle variables, then the integrals depend only on the action variables and the contact form takes on a similar form), as should be expected.
\end{remark}
\newpage
\chapter{Global aspects}
In this final chapter, we will take a look at some global properties regarding contact completely integrable systems. In particular, we will introduce contact toric $G$-manifolds (where $G$ is a certain torus), which is the contact version of toric systems in symplectic geometry. We will introduce the necessary concepts and end with a classification theorem of compact connected contact toric $G$-manifolds. These results are due to \cite{Lerman}. Morally, this chapter will coincide with chapter $3$ of the lecture notes on symplectic geometry.\newline

\noindent
Contrary to the fact that we will mostly consider torus actions, we will consider a general Lie group action wherever possible. For more information on Lie group (actions), we refer the reader to the lecture notes on symplectic geometry. Furthermore, we will only consider cooriented contact manifolds in this chapter. 
\section{Contact moment map}
Let $G$ be a Lie group acting on a (not necessarily contact) manifold $M$ via the diffeomorphisms $\varphi^G_g: M \rightarrow M$ with $g \in G$. We can always lift this action to the cotangent bundle $T^*M$ using the formula
\begin{equation*}
G\times T^*M \rightarrow T^*M: (g,(m,\beta_m)) \mapsto \left(\varphi^G_g(m),((\varphi^G_g)^{-1})^*(\beta_m)\right).
\end{equation*}
The following lemma will in particular be useful in the contact setting, since it discusses $1$-forms.
\newline \hfill
\newline
\noindent
\fbox{ \parbox{\textwidth}{
\begin{lemma}\label{Invariant Lemma}
Consider a manifold $M$ and a co-oriented codimension 1 distribution $\zeta$ on $M$ (with co-orientation $\zeta^\circ_+$). Let $G$ be a Lie group acting properly on $M$ such that $\zeta$ and its coorientation are preserved (in particular, the lifted action of $G$ to $T^*M$ preserves $\zeta^\circ_+$). Then there exists a $G$-invariant $1$-form $\beta$ on $M$ such that $\ker \beta = \zeta$ and $\beta(M) \subset \zeta^\circ_+$.
\end{lemma}}}
\begin{proof}
See Lemma $2.6$ in \cite{Lerman}. The proof is based on an averaging argument which is, depending on which exercises we decide to do in the regular course, outside of the scope of these notes.
\end{proof}
\begin{remark}\noindent
If $G$ acts properly on a contact manifold $(M,\xi = \ker \alpha)$ while preserving $\xi$ and its coorientation, then the lemma above says that we can always assume that $\alpha$ is $G$-invariant. Since a torus action is always proper, this will be a strong tool for our study of torus actions.
\end{remark}\noindent
If a Lie group preserves a $1$-form while acting on a manifold, there is a natural notion of a momentum map associated to it.
\newline \hfill
\newline
\noindent
\fbox{ \parbox{\textwidth}{.
\begin{definition}
Let $G$ be a Lie group acting on a manifold $M$ (not necessarily contact) while preserving a $1$-form $\beta$. We define the \textit{$\beta$-moment map} $\Psi_\beta: M \rightarrow \mathfrak{g}^*$, where $\mathfrak{g}^*$ is the dual of the Lie algebra $\mathfrak{g}$ of $G$, as
\begin{equation*}
\forall x \in M, \forall X \in \mathfrak{g}: \langle \Psi_\beta(x),X\rangle = \beta_x(X^{\sharp}(x)),
\end{equation*}
where $X^\sharp$ is the infinitesimal generator of the action on $M$ induced by $X$ and where $\langle \cdot, \cdot \rangle$ is the dual pairing.
\end{definition}}}
\newline \hfill \newline
\noindent
Consider now a Lie group $G$ acting properly on a contact manifold $(M,\xi)$ while preserving $\xi$ and its coorientation. By the previous lemma, there exists a $G$-invariant contact form $\alpha$ which defines $\xi$ and its coorientation. At first glance, it would be reasonable to define the ``contact moment map" as the $\alpha$-moment map $\Psi_\alpha$. However, there are some serious problems if we would do this. \newline

\noindent
Indeed, we have started with an action that preserves a contact structure and a coorientation. Hence we want the contact moment map to depend only on $\xi$ and $\xi^\circ_+$ and not on $\alpha$. If $f$ is a $G$-invariant smooth function on $M$, then $e^f\alpha$ is a $G$-invariant contact form on $M$ which defines $\xi$ and $\xi^\circ_+$. However, $\Psi_{e^f\alpha} = e^f\Psi_\alpha$, so this notion depends heavily on the choice of $G$-invariant contact form, which we do not want.\newline

\noindent
Luckily, there is another construction which will give us a moment map that only depends on the contact structure and its coorientation.
\newline \hfill
\newline
\noindent
\fbox{ \parbox{\textwidth}{
\begin{definition}
Let $G$ be a Lie group acting on a contact manifold $(M,\xi)$ while preserving $\xi$ and its coorientation $\xi_+^\circ$. Consider the lift of this action to the cotangent bundle $T^*M$ (as described in the beginning of this section), which is naturally symplectic (consider the Liouville form on $T^*M$). This lifted action will be Hamiltonian with momentum map $\Phi: T^*M\rightarrow \mathfrak{g}^*$ given by
\begin{equation*}
\forall q \in M,\forall p_q \in T_q^*M, \forall X \in \mathfrak{g}: \langle \Phi(q,p_q),X\rangle = \langle p_q , X^\sharp(q)\rangle,
\end{equation*}
where $X^\sharp$ is the infinitesimal generator of the action on $M$ induced by $X$. Since the action preserves the coorientation, we can restrict $\Phi$ to $\xi^\circ_+$. The \textit{contact moment map} is defined as this restriction $\Psi = \Phi|_{\xi^\circ_+}: \xi^\circ_+ \rightarrow \mathfrak{g}^*$.
\end{definition}}}
\begin{remark}
Depending on the choice of exercises in the regular course, we may not have seen the Liouville form on $T^*M$. It is nothing more than a specific $1$-form $\lambda$ on the cotangent bundle. For this form, it will hold that $d\lambda$ is a symplectic form. The exact definition of $\lambda$ will not be important for these notes: we only need the explicit formula for the momentum map $\Phi$. We refer to the literature (or to the regular exercises) for more details.
\end{remark} \noindent
Both of these concepts are related by the following proposition, which can be checked by direct calculations.
\newline \hfill
\newline
\noindent
\fbox{ \parbox{\textwidth}{
\begin{proposition} \label{Relatie moment maps}
Consider a cooriented contact manifold $(M,\xi)$ with coorientation $\xi^\circ_+$. Let a Lie group $G$ act on $M$ preserving $\xi$ and $\xi^\circ_+$. Assume that there is a $G$-invariant contact form $\alpha$ on $M$ such that $\alpha(M) \subset \xi^\circ_+$ and $\ker\alpha = \xi$. Then the $\alpha$-moment map $\Psi_\alpha$ and the contact moment map $\Psi$ for this action are related by
\begin{equation*}
\Psi \circ \alpha = \Psi_\alpha.
\end{equation*}
\end{proposition}}} \newline \hfill \newline \noindent
The proposition above is one of the reasons why we can consider $\Psi$ as a ``universal" moment map. We conclude this section by stating the following properties of an $\alpha$-moment map.
\newline \hfill
\newline
\noindent
\fbox{ \parbox{\textwidth}{
\begin{proposition}
Let $G$ be a Lie group which acts on a cooriented contact manifold $(M,\xi)$ preserving $\xi$ and its coorientation. Denote this action by the maps $\varphi^G_g$ for $g\in G$. Assume that $\alpha$ is a $G$-invariant contact form on $M$, defining $\xi$ and its coorientation. Let $\Psi_\alpha$ be the $\alpha$-moment map.
\begin{itemize}
\item[a)] It holds that $\Psi_\alpha$ is equivariant with respect to $\varphi^G$ and the coadjoint action, i.e.
\begin{equation*}
\forall g \in G: \textrm{Ad}_g^* \circ \Psi_\alpha = \Psi_\alpha \circ \varphi^G_g.
\end{equation*}
\item[b)] Consider the orbit $\mathcal{O}(m)$ of $\varphi^G$ through $m \in M$. For every $m \in \Psi_\alpha^{-1}(0)$, it holds that $T_m \mathcal{O}(m)$ is an isotropic subspace of $(\textrm{Ker}( \alpha_m), (d\alpha)_m)$.
\end{itemize}
\end{proposition}}}
\section{Contact toric manifolds}
Now, let us consider torus actions. Since a torus is compact, every torus action will be proper.
\newline \hfill
\newline
\noindent
\fbox{ \parbox{\textwidth}{
\begin{definition}
Let $(M,\xi)$ be a cooriented contact manifold and let $G$ be a torus. An action of $G$ on $(M,\xi)$ is called \textit{completely integrable} if
\begin{itemize}
\item The action is effective.
\item The action preserves $\xi$ and its coorientation.
\item $2\dim G = \dim M + 1$.
\end{itemize}
A cooriented contact manifold endowed with a completely integrable action of a torus $G$ is called a \textit{contact toric $G$-manifold}.
\end{definition}}}\newline \noindent \newline \noindent
The following lemma makes it possible to represent a contact toric $G$-manifold as a certain triple.
\newline \hfill
\newline
\noindent
\fbox{ \parbox{\textwidth}{
\begin{lemma}\label{Zero Image}
Consider a completely integrable action of a torus $G$ on a cooriented contact manifold $(M,\xi)$. Then $0$ is not in the image of the contact moment map $\Psi$.
\end{lemma}}}
\begin{proof}
See Lemma 2.12 in \cite{Lerman}, which makes use of the slice representation.
\end{proof}\noindent
Now, consider a cooriented contact manifold $(M,\xi)$ with a completely integrable action of a torus $G$. We can fix an inner product on the Lie algebra $\mathfrak{g}$, which induces a norm on the dual Lie algebra $\mathfrak{g}^*$. By lemma \ref{Invariant Lemma}, there exists a $G$-invariant contact form $\alpha'$ on $M$ such that $\ker \alpha' = \xi$ and $\alpha'(M) \subset \xi^\circ_+$.  By proposition \ref{Relatie moment maps} and lemma \ref{Zero Image}, it holds that $0$ is not in the image of $\Psi_{\alpha'}$. Hence we can define a new $G$-invariant contact form $\alpha$ on $M$ as
\begin{equation*}
\forall x \in M: \alpha_x = \frac{1}{||\Psi_{\alpha'}(x)||}\alpha_x',
\end{equation*}
such that it still holds that $\ker \alpha = \xi$ and $\alpha(M) \subset \xi^\circ_+$, but we also have that $||\Psi_\alpha(x)|| = 1$ for every $x \in M$. This special contact form $\alpha$ will be called the \textit{normalized contact form} for the contact toric $G$-manifold.
\newline
\noindent
By exercise \ref{Oef triple} from this chapter, a contact toric $G$-manifold is completely determined by the triple $(M,\alpha,\Psi_\alpha)$ (with $\alpha$ the normalized contact form). Hence from now on, we will denote a contact toric $G$-manifold as such a triple.\newline

\noindent Finally, we would like to state a classification theorem of contact toric $G$-manifolds. Before we can do this, we need to establish which contact toric $G$-manifolds are considered essentially the same.
\newline \hfill
\newline
\noindent
\fbox{ \parbox{\textwidth}{
\begin{definition} An \textit{isomorphism between two contact toric $G$-manifolds} $(M,\alpha,\Psi_\alpha)$ and $(M',\alpha',\Psi_{\alpha'})$ is a map $\varphi: M \rightarrow M'$ such that
\begin{itemize}
\item $\varphi$ is a contactomorphism.
\item $\varphi$ preserves the coorientation.
\item $\varphi$ is $G$-equivariant.
\end{itemize}
If such a map exists, we call  $(M,\alpha,\Psi_\alpha)$ and $(M',\alpha',\Psi_{\alpha'})$ \textit{isomorphic}.
\end{definition}}}\newline \hfill \newline \noindent
\begin{remark}\label{Remark isomorfisme}
The normalization conditions of $\alpha$ and $\alpha'$ forces every such contactomorphism $\varphi$ to be in fact a strict contactomorphism.
\end{remark}\noindent
The following geometric notions will turn up in the classification theorem:
\newline

\noindent
\fbox{ \parbox{\textwidth}{
\begin{definition}
Consider a cooriented contact manifold $(M,\xi)$ with the action of a Lie group $G$. Assume that this action preserves the contact structure and its coorientation. Let $\Psi$ be the contact moment map associated to this action. The \textit{moment cone} is defined as the set
\begin{equation*}
C(\Psi) = \Psi(\xi_+^\circ) \cup \{0\}.
\end{equation*}
In particular, if $\alpha$ is a $G$-invariant contact form which defines $\xi$ and its coorientation, we find that
\begin{equation*}
C(\Psi) = \{t\Psi_\alpha(x)|x \in M, t \geq 0\}.
\end{equation*}
\end{definition}
}}
\newline \hfill
\newline \hfill \newline
\noindent
\fbox{ \parbox{\textwidth}{

\noindent
\begin{definition}
Let $G$ be a torus and let $\mathfrak{g}^*$ be the dual of its Lie algebra. We say that a subset $C\subset \mathfrak{g}^*$ is a \textit{rational polyhedral cone} if there exists a finite set of vectors $\{v_i\}$ in the integral lattice $\mathbb{Z}_G$ such that $$C = \bigcap_{i} \{\eta \in \mathfrak{g}^* | \langle \eta,v_i\rangle \geq 0\}.$$
We can assume wlog that this set is minimal (if you remove a vector $v_i$, the intersection is strictly bigger than $C$) and that the vectors are primitive ($tv_i \notin \mathbb{Z}_G$ if $0 < t < 1$, so a shorter vector is not in $\mathbb{Z}_G$).\\
Let $C$ be a rational polyhedral cone, defined by the vectors $\{v_i\}$ in $\mathbb{Z}_G$. We say that $C$ is a \textit{good} rational polyhedral cone if
\begin{itemize}
\item The interior of $C$ is non-empty.
\item Consider the linear span of a codimension $k$ face (for $0<k<\dim G)$ of $C$. The annihilator of this space (which lives in $\mathfrak{g}$) is the Lie algebra of a subtorus $H$ of the original torus $G$.
\item The normal vectors of the face described above (i.e. the $k$ different vectors $v_{i_j}$ from $\{v_i\}$ such that $\langle \eta,v_{i_j}\rangle = 0$ for every $\eta$ in the face) form a $\mathbb{Z}$-basis of the integral lattice $\mathbb{Z}_H$ of $H$.
\end{itemize}
\end{definition}}}
\newpage \noindent
\begin{example}
Consider $G = \mathbb{T}^2$ and thus $\mathfrak{g}^* = \mathbb{R}^2$ and $\mathbb{Z}_G = \mathbb{Z}^2$. We define $$v_1 = \begin{bmatrix}
0 \\ -1
\end{bmatrix} \quad \textrm{and} \quad v_2 =  \begin{bmatrix}
1 \\ 1
\end{bmatrix}.$$
Then we have that the set $$ C = \{\eta \in \mathbb{R}^2 | \langle \eta,v_1\rangle \geq 0, \langle \eta,v_2\rangle \geq 0\}$$ is a rational polyhedral cone. It is a direct check to see that $C$ is in fact a good rational polyhedral cone.
\end{example}\noindent
The following special class of manifolds will appear in the classification.
\newline \hfill
\newline
\noindent
\fbox{ \parbox{\textwidth}{
\begin{definition}
Consider $\mathbb{S}^3$ as a submanifold of $\mathbb{C}^2$. Let $p,q \in \mathbb{Z}$ be coprime integers. Consider the (free) action of the group $\mathbb{Z}/p\mathbb{Z}$ on $\mathbb{S}^3$, generated by the rotation
\begin{equation*}
\rho: \mathbb{S}^3 \rightarrow \mathbb{S}^3: (z_1,z_2) \mapsto (e^{\frac{2\pi i}{p}}z_1, e^{\frac{2\pi iq}{p}}z_2).
\end{equation*}
The resulting orbit space $\mathbb{S}^3/(\mathbb{Z}/p\mathbb{Z})$ is called a \textit{lens space} and is denoted by $L(p,q)$ (notice that the lens space is determined by the integers $(p,q)$).
\end{definition}}}
\begin{example}
Some examples of lens spaces are
\begin{itemize}
\item $L(0,1) = \mathbb{S}^1 \times \mathbb{S}^2$
\item $L(1,0) = \mathbb{S}^3$
\item $L(2,1) = \mathbb{RP}^3$
\end{itemize}
As an example, let us prove that $L(2,1) = \mathbb{RP}^3$. The $\mathbb{Z}/2\mathbb{Z}$-action on $\mathbb{S}^3$ (considered as submanifold of $\mathbb{C}^2$) is
\begin{equation*}
\mathbb{Z}/2\mathbb{Z} \times \mathbb{S}^3 \rightarrow \mathbb{S}^3: (m,(z_1,z_2)) \mapsto (e^{\frac{2\pi m i}{2}}z_1, e^{\frac{2\pi m i\cdot 1}{2}}z_2) = (e^{\pi m i}z_1, e^{\pi m i}z_2).
\end{equation*}
Since $m$ is either $0$ or $1$, the orbit of a point $(z_1,z_2)$ under this action is $\{(z_1,z_2),(-z_1,-z_2)\}$. Consequently, the orbit space $\mathbb{S}^3/(\mathbb{Z}/2\mathbb{Z})$ (which is exactly the lens space $L(2,1)$) is the $3$-sphere $\mathbb{S}^3$ where every point $(z_1,z_2)$ is identified with its antipodal point $(-z_1,-z_2)$. This is a well-known model of the projective space $\mathbb{RP}^3$. We conclude that $L(2,1) = \mathbb{RP}^3$.
\end{example}
\noindent We can now finally state the classification result.
\newpage
\noindent
\fbox{ \parbox{\textwidth}{
\begin{theorem}
(\textbf{Lerman classification}) Compact connected contact toric $G$-manifolds (abbreviated in this theorem as c.c.c.t. $G$-manifolds) $(M,\alpha,\Psi_\alpha)$ for a torus $G$ can be classified in the following way:
\begin{itemize}
\item[i)] Let $\dim M = 3$ and assume the $G$-action (with $G = \mathbb{T}^2$) is free. Then $M$ is diffeomorphic to $\mathbb{T}^3 = \mathbb{S}^1 \times \mathbb{T}^2$. The contact form $\alpha$ is of the form 
\begin{equation*}
\alpha = \cos (nt) d\theta_1 + \sin(nt) d\theta_2,
\end{equation*}
where $(t,\theta_1,\theta_2) \in \mathbb{S}^1 \times \mathbb{T}^2$ and $n \in \mathbb{N}\setminus \{0\}$.
\item[ii)] Let $\dim M = 3$ and assume the $G$-action (with $G = \mathbb{T}^2$) is not free. Then $M$ is diffeomorphic to a lens space. Furthermore, as a c.c.c.t. $G$-manifold, the triple $(M,\alpha,\Psi_\alpha)$ is classified by two rational numbers $r,q$ such that $ 0\leq r < 1$ and $r < q$.
\item[iii)] Let $\dim M > 3$ and assume the $G$-action is free. Let $d = \dim G -1$. It holds that $M$ is a principal $G$-bundle over $\mathbb{S}^d$. Furthermore, every principal $G$-bundle over $\mathbb{S}^d$ has a unique $G$-invariant contact structure which makes this bundle a c.c.c.t. $G$-manifold.
\item[iv)] Let $\dim M > 3$ and assume the $G$-action is not free. In this case, the moment cone of $(M,\alpha,\Psi_\alpha)$ is a good cone.\\
On the other hand, if $C \subset \mathfrak{g}^*$ is a good cone, there exists a unique c.c.c.t. $G$-manifold $(M,\alpha,\Psi_\alpha)$ with moment cone $C$.
\end{itemize}
\end{theorem}}}
\begin{proof}
See theorem $2.18$ in \cite{Lerman}, which is the main result of that paper. In very broad terms, he first classifies neighbourhoods of so called pre-isotropic embeddings. Then, he makes this classification global by considering the Cech cohomology of orbit spaces with coefficients in a certain abelian sheaf.
\end{proof}
\noindent
This theorem has $4$ different subcases, but the essence is that not a lot can happen. In the first case, the manifold has to be a torus and there are countably many contact forms which make it a c.c.c.t. $G$-manifold (parametrized by $n$). In the second case, the manifold has to be a lens space and there are again only countably many options for this space to be a c.c.c.t. $G$-manifold (parametrized by $r$ and $q$). In the third case, there is a notion from differential geometry that we have not discussed\footnote{I refer to the course Selected topics in differential geometry for more details.}, but it can be shown that for $\dim M > 5$, there is only $1$ c.c.c.t. $G$-manifold and for $\dim M = 5$, there are countably many (namely $\mathbb{Z}^3$). Finally, the fourth case comes down to performing the Delzant construction on the given cone and restricting to a hypersurface of contact type.\newline

\noindent
The discussion above shows that the world of contact toric systems is rather rigid. Hence a natural question to ask is if there exists a class of (contact) systems which are less rigid than contact toric systems, but which are still manageable to classify. In particular, we would like to know if there exists something like \textit{contact semitoric systems}. This notion has not yet been defined (let alone classified). However, the introduction and classification of such systems forms an active part of the research performed by the Analysis group at the University of Antwerp.
\section{Lerman construction}
Let us now focus on the final case of Lerman's classification. In this section, we will see how one can construct a compact connected contact toric $G$-manifold from a good cone (given that the dimension of the underlying manifold is strictly bigger than 3). The general idea is that we first perform Delzant's construction (compare this to lecture notes on symplectic geometry), after which we restrict to a hypersurface of contact type (see proposition \ref{contactisation}). We will focus on the results from \cite{Lerman}.\\ \hfill \\
Consider a torus $G$ with dimension $n \geq 3$ with Lie algebra $\mathfrak{g}$. Let $\left\{v_j\right\}_{1\leq j \leq d}$ be a finite set of vectors in the integral lattice $\mathbb{Z}_G$ such that 
\begin{equation*}
C = \bigcap_{j=1}^d \{\eta \in \mathfrak{g}^* | \langle \eta,v_j\rangle \geq 0\}
\end{equation*}
is a good rational polyhedral cone. We may assume that the set $\left\{v_j\right\}_{1\leq j \leq d}$ is minimal and that its elements are primitive.\\
Our first goal is to construct a symplectic manifold $(W,\omega)$ together with a Hamiltonian action of the torus $G$ (in the sense of the definition given in the lecture notes on symplectic geometry) such that $C$ is the image of the associated momentum map. To achieve this goal, we will follow Delzant's contruction.\\
Consider the standard basis $\left\{e_j\right\}_{1 \leq j \leq d}$ of $\mathbb{R}^d$. We define a linear map $\tau$ by sending $e_k$ to $v_k$, i.e.
$$\tau: \mathbb{R}^d \rightarrow \mathfrak{g} \simeq \mathbb{R}^n: \sum_{j=1}^d a_j e_j \mapsto \sum_{j=1}^d a_j v_j.$$
If we would have performed a linear extension over $\mathbb{Z}$ instead, we would get a well-defined map $\tau_\mathbb{Z}: \mathbb{Z}^d \rightarrow \mathbb{Z}_G \simeq \mathbb{Z}^n$ (since $\tau(\mathbb{Z}^d) \subset \mathbb{Z}_G$). Hence the map $\tau$ descends to the quotient:
$$\widetilde{\tau}: \mathbb{T}^d \simeq\mathbb{R}^d/\mathbb{Z}^d \rightarrow \mathfrak{g}/\mathbb{Z}_G \simeq G \simeq \mathbb{T}^n: \widetilde{a} \mapsto \widetilde{\tau(a)},$$
where we use the following convention: if $u \in \mathbb{R}^m$ then $\widetilde{u} \in \mathbb{R}^m/\mathbb{Z}^m \simeq \mathbb{T}^m$ is the image of $u$ under the natural projection map. We will denote the kernel of $\widetilde{\tau}$ as $N$, so $$N = \left\{ \widetilde{a} \in \mathbb{T}^d| \sum_{j=1}^d a_jv_j \in \mathbb{Z}_G \right\}.$$
One can directly check that $N$ is a compact, abelian subgroup of $\mathbb{T}^d$, so the connected component of $N$ that contains the neutral element is a torus itself. Denote the Lie algebra of $N$ as $\mathfrak{n}$, which is equal to the kernel of $\tau$. Furthermore, define the inclusion map $i: \mathfrak{n} \rightarrow \mathbb{R}^d$, which we will need later.\\
Just as in the lecture notes on symplectic geometry, we now consider the standard $\mathbb{T}^d$-action $\varphi$ on $\left(\mathbb{C}^d,\omega_0 = \frac{i}{2}\sum_{k=1}^d dz_k \wedge d\overline{z_k}\right)$, defined as
\begin{equation*}
\varphi: \mathbb{T}^d \times \mathbb{C}^d \rightarrow \mathbb{C}^d: \left((\widetilde{t_1},\cdots,\widetilde{t_d}),(z_1\cdots,z_d)\right) \rightarrow \left(e^{2\pi it_1}z_1,\cdots,e^{2\pi it_d}z_d \right).
\end{equation*}
One can prove as an exercise that this action is Hamiltonian with momentum map 
\begin{equation*}
h: \mathbb{C}^d \rightarrow \textrm{Lie}\left(\mathbb{T}^d \right)^* \simeq \mathbb{R}^d: (z_1,\cdots,z_d) \mapsto \pi \sum_{j=1}^d |z_j|^2e_j^*,
\end{equation*}
where $\left\{e_j^*\right\}_{1 \leq j \leq d}$ is the basis of $\textrm{Lie}\left(\mathbb{T}^d \right)^*$ dual to $\left\{e_j\right\}_{1 \leq j \leq d}$. Now, define the map
$$\widetilde{h} = i^* \circ h: \mathbb{C}^d \rightarrow \mathfrak{n}^* \simeq \mathbb{R}^{d-n}.$$
In order to take the Marsden-Weinstein quotient, we need to check following technical properties.\\ \hfill \\
\noindent
\fbox{ \parbox{\textwidth}{
\begin{lemma}
The action of the Lie group $N$ on $\mathbb{C}^d$, induced by $\varphi$, leaves the level set $\widetilde{h}^{-1}(0)$ invariant.
\end{lemma}}}
\begin{proof}\noindent
We will show that the action of $\mathbb{T}^d$ will leave $\widetilde{h}^{-1}(0)$ invariant, so it will surely hold for $N$.  Let $\widetilde{t} = (\widetilde{t_1},\cdots,\widetilde{t_d}) \in \mathbb{T}^d$ and $z = (z_1\cdots,z_d) \in \widetilde{h}^{-1}(0)$. It holds that
\begin{eqnarray*}
\widetilde{h}\left(\varphi(\widetilde{t},z)\right)&=& i^*\left(h\left(\varphi(\widetilde{t},z)\right) \right)\\
&=& i^*\left(h\left(e^{2\pi it_1}z_1,\cdots,e^{2\pi it_d}z_d\right) \right)\\
&=& i^*\left(\pi \sum_{j=1}^d |e^{2\pi it_j}z_j|^2e_j^*\right)\\
&=& i^*\left(\pi \sum_{j=1}^d |z_j|^2e_j^*\right)\\
&=& i^*\left(h(z)\right)\\
&=& \widetilde{h}(z)\\
&=& 0
\end{eqnarray*}
We conclude that $\varphi(\widetilde{t},z)\in \widetilde{h}^{-1}(0)$.
\end{proof}\noindent
Consequently, we have an action of $N$ on $\widetilde{h}^{-1}(0)$.\\ \hfill \\
\noindent
\fbox{ \parbox{\textwidth}{
\begin{lemma}
The $N$-action on $\widetilde{h}^{-1}(0)$ is proper.
\end{lemma}}}
\begin{proof}
This is clear, since $N$ is a compact Lie group and every action of a compact Lie group on a smooth manifold is proper.
\end{proof}
\noindent
\fbox{ \parbox{\textwidth}{\label{technisch lemma Lerman constructie}
\begin{lemma}
$\widetilde{h}^{-1}(0) = h^{-1}(\tau^*(C))$.
\end{lemma}}}
\begin{proof}
Since the sequence
\begin{equation*}
\begin{tikzcd}              
    0 \arrow{r} & \mathfrak{n} \arrow[r,"i"]& \mathbb{R}^d  \arrow[r,"\tau"]       & \mathfrak{g}           \arrow{r} & 0
\end{tikzcd}
\end{equation*}
is a short exact sequence (where we use that $C$ is a good rational polyhedral cone), so is the ``dual" sequence
\begin{equation*}
\begin{tikzcd}              
    0 \arrow{r} & \mathfrak{g}^* \arrow[r,"\tau^*"]& (\mathbb{R}^d)^*                 \arrow[r,"i^*"]       & \mathfrak{n}^*           \arrow{r} & 0.
\end{tikzcd}
\end{equation*}
Hence $(i^*)^{-1}(0) = \tau^*(\mathfrak{g}^*)$. Consequently,
$$\widetilde{h}^{-1}(0) = h^{-1}\left((i^*)^{-1}(0)\cap h(\mathbb{C}^d)\right) = h^{-1}\left(\tau^*(\mathfrak{g}^*)\cap h(\mathbb{C}^d)\right).$$
Furthermore, since $$h\left(\mathbb{C}^d\right) = \left\{\eta \in (\mathbb{R}^d)^* | \langle\eta,e_j\rangle \geq 0 \textrm{ for } 1\leq j \leq d \right\},$$
we can find that
\begin{eqnarray*}
\tau^*(\mathfrak{g}^*)\cap h(\mathbb{C}^d) &=& \left\{\tau^*(\eta)| \eta \in \mathfrak{g}^* \textrm{ and } \langle\tau^*(\eta) ,e_j\rangle \geq 0 \textrm{ for } 1\leq j \leq d \right\}\\
&=& \left\{\tau^*(\eta)| \eta \in \mathfrak{g}^* \textrm{ and } \langle\eta ,\tau(e_j)\rangle \geq 0 \textrm{ for } 1\leq j \leq d \right\}\\
&=& \left\{\tau^*(\eta)| \eta \in \mathfrak{g}^* \textrm{ and } \langle\eta ,v_j\rangle \geq 0 \textrm{ for } 1\leq j \leq d \right\}\\
&=& \left\{\tau^*(\eta)| \eta \in C \right\}\\
&=& \tau^*\left( C \right).
\end{eqnarray*}
We can finally conclude that $\widetilde{h}^{-1}(0) = h^{-1}\left(\tau^*\left( C \right)\right)$.
\end{proof}\noindent
This result can be used to prove that the action of $N$ on $\widetilde{h}^{-1}(0)\setminus \left\{0\right\}$ is free.\\ \hfill \\
\noindent
\fbox{ \parbox{\textwidth}{\label{Lerman construction triviale stabilisator}
\begin{lemma}
Let $z \in \widetilde{h}^{-1}(0)\setminus \left\{0\right\}$. It holds that the isotropy group $N_z$ of $z$ with respect to the action of $N$ on $\widetilde{h}^{-1}(0)\setminus \left\{0\right\}$ is trivial, i.e. $N_z = \left\{0\right\}$.
\end{lemma}}}
\begin{proof}
Consider the isotropy group $(\mathbb{T}^d)_z$ of $z$ with respect to the $\mathbb{T}^d$-action on $\widetilde{h}^{-1}(0)\setminus \left\{0\right\}$. We have that
\begin{eqnarray*}
N_z &=& N \cap (\mathbb{T}^d)_z \\
&=& \left\{ \widetilde{a} \in \mathbb{T}^d| \sum_{j=1}^d a_jv_j \in \mathbb{Z}_G \right\} \cap \left\{ \widetilde{a} \in \mathbb{T}^d| a_j \in \mathbb{Z} \textrm{ if } j\in \left\{1,\cdots,d\right\} \textrm{ such that } z_j \neq 0 \right\}\\
&=& \left\{ \widetilde{a} \in \mathbb{T}^d| \sum_{j\in J_z} a_jv_j \in \mathbb{Z}_G \textrm{ and } a_k \in \mathbb{Z} \textrm{ if } k \notin J_z\right\},
\end{eqnarray*}
where we define $J_z$ as $$J_z = \left\{j \in \left\{1,\cdots, d\right\}| z_j = 0 \right\}.$$
However, using lemma \ref{technisch lemma Lerman constructie}, there exists an $\eta \in C$ such that $h(z) = \tau^*(\eta)$. This means that $z_j = 0$ if and only if $$0 = \pi|z_j|^2 =  \langle h(z),e_j \rangle = \langle \tau^*(\eta),e_j \rangle = \langle \eta,\tau(e_j) \rangle =  \langle \eta,v_j \rangle.$$
We can thus rewrite $J_z$ as $$J_z = \left\{j \in \left\{1,\cdots, d\right\}|\langle \eta,v_j \rangle = 0 \right\}.$$
Using the fact that $C$ is a good cone, we have that $\left\{v_j \right\}_{j \in J_z}$ is a $\mathbb{Z}$-basis of $$\left\{ \sum_{j\in J_z}a_jv_j | a_j \in \mathbb{R}\right\}\cap \mathbb{Z}_G.$$
Consequently, if $ \sum_{j\in J_z}a_jv_j \in \mathbb{Z}_G$ (in particular, if $\widetilde{a} \in N_z$), then $a_j \in \mathbb{Z}$ for every $j \in J_z$. Putting everything together, we get
$$ N_Z = \left\{ \widetilde{a} \in \mathbb{T}^d| a \in \mathbb{Z}^N\right\} = \left\{0\right\},$$
so $N_z$ is trivial.
\end{proof}
\noindent
\fbox{ \parbox{\textwidth}{
\begin{corollary}
The action of $N$ on $\widetilde{h}^{-1}(0)\setminus \left\{0\right\}$ is free.
\end{corollary}}}
\begin{proof}
This follows directly from lemma \ref{Lerman construction triviale stabilisator}.
\end{proof}\noindent
Since we have a smooth, proper and free action of $N$ on $\widetilde{h}^{-1}(0)\setminus \left\{0\right\}$, we can take the Marsden-Weinstein quotient to get a symplectic manifold $(M_{red,0},\omega)$, where $$M_{red,0} = \left(\widetilde{h}^{-1}(0)\setminus \left\{0\right\} \right)/N$$ and where $\omega$ is implicitly defined as
$$q^*\omega = j^*\omega_0,$$
with the natural inclusion map $j: \widetilde{h}^{-1}(0)\setminus \left\{0\right\} \rightarrow \mathbb{C}^d$ and the natural projection map $q:\widetilde{h}^{-1}(0)\setminus \left\{0\right\} \rightarrow \left(\widetilde{h}^{-1}(0)\setminus \left\{0\right\}\right)/N$.\\
We will now construct a momentum map for the induced $\mathbb{T}^d/N \simeq G$-action on $(M_{red,0},\omega)$, just as we have done in the symplectic case. Let $\sigma: \mathfrak{g} \rightarrow \mathbb{R}^d$ be a right inverse of $\tau$. Then we can define a map $$\widetilde{F} = \sigma^* \circ h \circ j: \widetilde{h}^{-1}(0)\setminus \left\{0\right\} \rightarrow \mathfrak{g}^* \simeq \mathbb{R}^n.$$
As we have seen in the Delzant construction, this map is $N$-invariant, so it descends to the quotient:
$$F: M_{red,0} = \left(\widetilde{h}^{-1}(0)\setminus \left\{0\right\}\right)/N \rightarrow \mathfrak{g}^* \simeq \mathbb{R}^n. $$
This map is the momentum map for the induced $G$-action on $(M_{red,0},\omega).$ The image of this map will be the cone $C$ without $0$. \\ \hfill \\
\noindent
\fbox{ \parbox{\textwidth}{
\begin{proposition}
$F(M_{red,0}) = C\setminus \{0\}$.
\end{proposition}}}
\begin{proof}
We will first compute the image of $\widetilde{F}$. Using lemma \ref{technisch lemma Lerman constructie}, we have that
\begin{eqnarray*}
\widetilde{F}\left(\widetilde{h}^{-1}(0)\setminus \left\{0\right\}\right) &=& \left(\sigma^* \circ h \circ j\right)\left(\widetilde{h}^{-1}(0)\setminus \left\{0\right\}\right)\\
&=& \sigma^*\left(h\left(\widetilde{h}^{-1}(0)\setminus \left\{0\right\}\right)\right)\\
&=& \sigma^*\left(h\left(h^{-1}\left(\tau^*(C\setminus \{0\})\right)\right)\right)\\
&=& \sigma^*\left(\tau^*(C\setminus \{0\})\right)\\
&=& \left(\tau \circ \sigma \right)^*(C\setminus \{0\})\\
&=& C\setminus \{0\},
\end{eqnarray*}
where we have used that $\sigma$ is a right inverse of $\tau$. Since $\widetilde{F}$ is an $N$-invariant map, it holds that $F(M_{red,0}) = C\setminus \{0\}$.
\end{proof}\noindent
In the second part of Lerman's construction, we will restrict to a certain hypersurface of contact type. We define an $\mathbb{R}$-action on $\mathbb{C}^d$:
$$\psi: \mathbb{R} \times \mathbb{C}^d \rightarrow \mathbb{C}^d: (t,z) \mapsto e^{\frac{t}{2}}z.$$
\noindent
\fbox{ \parbox{\textwidth}{
\begin{lemma}\label{Lerman construction Liouville VF voor quotient}
Let $\psi_t: \mathbb{C}^d \rightarrow \mathbb{C}^d: z \mapsto \psi(t,z)$. It holds that $\psi_t^*\omega_0 = e^t \omega_0$. 
\end{lemma}}}
\begin{proof}
This can be shown with a direct calculation:
\begin{eqnarray*}
\psi_t^* \omega_0 &=& \psi_t^* \left(\frac{i}{2}\sum_{k=1}^d dz_k \wedge d\overline{z}_k\right)\\
&=&\frac{i}{2}\sum_{k=1}^d d\left(z_k \circ \psi_t\right) \wedge d\left(\overline{z}_k \circ \psi_t\right)\\
&=&\frac{i}{2}\sum_{k=1}^d d\left(e^{\frac{t}{2}}z_k\right) \wedge d\left(\overline{e^{\frac{t}{2}}z_k}\right)\\
&=&e^t\frac{i}{2}\sum_{k=1}^d dz_k \wedge d\overline{z}_k\\
&=& e^t\omega_0
\end{eqnarray*}
Hence $\psi_t^*\omega_0 = e^t \omega_0$.
\end{proof}\noindent
By exercise \ref{oef implicatie Liouville} from chapter \ref{Chapter Contact geometry}, we know that the vector field $A$, defined as $$\forall z \in \mathbb{C}^d: A(z) = \frac{\textrm{d}}{\textrm{d}t}\Big|_{t=0} \psi(t,z)$$ is a Liouville vector field on $(\mathbb{C}^d,\omega_0)$. Furthermore, it is clear that $A$ is transverse to the hypersurface
$$\mathbb{S}^{2d-1} = \left\{z \in \mathbb{C}^d | \sum_{j=1}^d |z_j|^2 = 1\right\},$$
since its flow is a rescaling in the radial direction. Consequently, $\mathbb{S}^{2d-1}$ is a hypersurface of contact type of the symplectic manifold $(\mathbb{C}^d,\omega_0)$.\\
This hypersurface can be used to find a hypersurface of contact type of $(M_{red,0},\omega)$. First, we need to show that the action $\psi$ descends to the quotient.
\\ \hfill \\
\noindent
\fbox{ \parbox{\textwidth}{
\begin{lemma}\label{Lerma construction lemma actie daalt neer}
$\psi$ induces an action on $\widetilde{h}^{-1}(0)\setminus \left\{0\right\}$.
\end{lemma}}}
\begin{proof}
Let $z \in \widetilde{h}^{-1}(0)\setminus \left\{0\right\}$ and let $t \in \mathbb{R}$. We have that
\begin{eqnarray*}
\widetilde{h}(\psi(t,z)) &=& \widetilde{h}\left(e^{\frac{t}{2}}z\right)\\
&=& i^*\left(h\left(e^{\frac{t}{2}}z\right)\right)\\
&=& i^*\left(\pi \sum_{j=1}^d |e^{\frac{t}{2}}z_j|^2e_j^*\right)\\
&=& i^*\left(e^t\pi \sum_{j=1}^d |z_j|^2e_j^*\right)\\
&=& i^*\left(e^th(z)\right).
\end{eqnarray*}
Let $X \in \mathbb{R}^d$. We can show that we can take $e^t$ outside of $i^*$:
$$\langle i^*\left(e^th(z)\right), X \rangle = \langle e^th(z), i(X) \rangle = e^t \langle h(z), i(X) \rangle = \langle e^ti^*(h(z)), X \rangle,$$
so $ i^*\left(e^th(z)\right) = e^ti^*(h(z))$ by the non-degeneracy of the dual pairing. So  $$\widetilde{h}(\psi(t,z)) = e^ti^*(h(z)) = e^t\widetilde{h}(z) = 0.$$
This shows that $\psi$ induces an action on $\widetilde{h}^{-1}(0)\setminus \left\{0\right\}$.
\end{proof}
\noindent
\fbox{ \parbox{\textwidth}{
\begin{corollary}
The action $\psi$ descends to the quotient, i.e. we have an action $$\widetilde{\psi}: \mathbb{R} \times M_{red,0} \rightarrow M_{red,0}$$
that satisfies $q \circ \psi_t = \widetilde{\psi}_t \circ q$ for every $t \in \mathbb{R}$.
\end{corollary}}}
\begin{proof}
This follows immediately from lemma \ref{Lerma construction lemma actie daalt neer} and the fact that the actions $\varphi$ and $\psi$ commute.
\end{proof}\noindent
This action will now generate a Liouville vector field on $(M_{red,0},\omega)$.\\ \hfill \\
\noindent
\fbox{ \parbox{\textwidth}{
\begin{proposition}
Consider the vector field $Y$ on $M_{red,0}$, defined as
$$\forall z \in M_{red,0}: Y(z) = \frac{\textrm{d}}{\textrm{d}t}\Big|_{t=0} \widetilde{\psi}(t,z).$$
Then $Y$ is a Liouville vector field on $M_{red,0}$.
\end{proposition}}}
\begin{proof}
We will show that $\widetilde{\psi}_t^*\omega = e^t \omega$, after which the result follows from exercise \ref{oef implicatie Liouville} from chapter \ref{Chapter Contact geometry}. We make the following computation:
\begin{eqnarray*}
q^* \left(\frac{1}{e^t}\widetilde{\psi}_t^*\omega\right) &=& q^* \left(\widetilde{\psi}_t^*\left(\frac{1}{e^t}\omega\right)\right)\\
&=& \left(\widetilde{\psi}_t \circ q \right)^*\left(\frac{1}{e^t}\omega\right)\\
&=& \left(q \circ \psi_t \right)^*\left(\frac{1}{e^t}\omega\right)\\
&=& \psi_t^*\left( q^*\left(\frac{1}{e^t}\omega\right)\right)\\
&=&  \psi_t^*\left(\frac{1}{e^t} q^*\omega\right)\\
&=& \psi_t^*\left(\frac{1}{e^t} j^*\omega_0\right)\\
&=& \frac{1}{e^t} (j\circ \psi_t)^*\omega_0
\end{eqnarray*}
Recall that $j$ is a natural inclusion map, hence (by abuse of notation) $$j \circ \psi_t = \psi_t \circ j.$$
If we now apply lemma \ref{Lerman construction Liouville VF voor quotient}, we find that $$q^* \left(\frac{1}{e^t}\widetilde{\psi}_t^*\omega\right) = \frac{1}{e^t} (\psi_t\circ j)^*\omega_0 = \frac{1}{e^t} j^*\left(\psi_t^*\omega_0\right) = \frac{1}{e^t} j^*\left(e^t\omega_0\right) = j^*\omega_0. $$
Since $\omega$ is implicitly defined as $q^*\omega = j^*\omega_0$, we have that $\frac{1}{e^t}\widetilde{\psi}_t^*\omega = \omega$, or equivalently $\widetilde{\psi}_t^*\omega = e^t\omega$. 
\end{proof}\noindent
Since $\mathbb{S}^{2d-1}$ is a $\mathbb{T}^d$-invariant hypersurface of $\mathbb{C}^d$ (with respect to the action $\varphi$), we have that $$M = \left(\widetilde{h}^{-1}(0) \cap \mathbb{S}^{2d-1} \right)/N$$
is a $G$-invariant hypersurface of $M_{red,0}$ (since $ G \simeq \mathbb{T}^d/N$). Furthermore, the Liouville vector field $Y$ is transverse to the hypersurface $M$ (since $A$ is transverse to $\mathbb{S}^{2d-1}$). Consequently, $M$ is a $G$-invariant hypersurface of contact type of $(M_{red,0},\omega)$, which can be endowed with the contact form $$\alpha = \iota_Y\omega.$$
Using exercise \ref{oef questionable} from this chapter and the fact that $F$ is the momentum map for the $G$-action on $(M_{red,0},\omega)$, we have that $$H = F|_M: M \rightarrow \mathfrak{g}^*$$
is the corresponding $\alpha$-moment map of the contact toric $G$-action on $(M,\ker\alpha)$. Furthermore, since $F(M_{red,0}) = C\setminus \{0\}$, we have that the moment cone of $(M,\alpha,H)$ is exactly $C$:
$$C(\Psi) = \{tH(x)|x \in M, t \geq 0\} = C,$$
where $\Psi$ is the contact moment map of the contact toric $G$-manifold $(M,\alpha,H)$. This is exactly what we wanted to construct.\newpage
\section{Exercises}
\begin{enumerate}

\item Let $G$ be a Lie group acting on a manifold $M$ while preserving a $1$-form $\beta$. Assume that $d\beta$ is a symplectic form on $M$. Show that $\Psi_\beta$ is a momentum map of the Hamiltonian action of $G$ on $(M,d\beta)$.

\item Consider a cooriented contact manifold $(M,\xi)$ with coorientation $\xi^\circ_+$. Let a Lie group $G$ act on $M$ while preserving $\xi$ and $\xi^\circ_+$. Denote the action by the maps $\varphi^G_g$ for $g\in G$. Assume that there exists a $G$-invariant contact form $\alpha$ on $M$ such that $\alpha(M) \subset \xi^\circ_+$ and $\ker\alpha = \xi$. Let $\Psi_\alpha$ be the $\alpha$-moment map.
\begin{enumerate}
\item Show that $\Psi_\alpha$ and the contact moment map $\Psi$ for this action are related by
\begin{equation*}
\Psi \circ \alpha = \Psi_\alpha.
\end{equation*}
\item Show that $\Psi_\alpha$ is equivariant with respect to $\varphi^G$ and the coadjoint action, i.e.
\begin{equation*}
\forall g \in G: \textrm{Ad}_g^* \circ \Psi_\alpha = \Psi_\alpha \circ \varphi^G_g.
\end{equation*}
\noindent
\textit{\underline{Hint:}} You can use the following property of the exponential map (without proof), with $\mathfrak{g}$ the Lie algebra of $G$:
\begin{equation*}
\forall X \in \mathfrak{g},\forall g \in G, \forall t \in \mathbb{R}: g\exp(tX)g^{-1} = \exp(t\textrm{Ad}(g)(X)).
\end{equation*}
\item Consider the orbit $\mathcal{O}(m)$ of $\varphi^G$ through $m \in M$. For every $m \in \Psi_\alpha^{-1}(0)$, it holds that $T_m \mathcal{O}(m)$ is an isotropic subspace of $(\textrm{Ker}( \alpha_m), (d\alpha)_m)$.
\end{enumerate}

\item \label{Oef triple} Consider a cooriented contact manifold $(M,\xi)$ with a completely integrable action of a torus $G$. Let $\alpha$ be the normalized contact form for this contact toric $G$-manifold. Show that one can reconstruct the entire action from the triple $(M,\alpha,\Psi_\alpha)$, with $\Psi_\alpha$ the $\alpha$-momentum map.

\item Prove the statement in remark \ref{Remark isomorfisme}. In particular, let $\varphi:M \rightarrow M'$ be an isomorphism between the contact toric $G$-manifolds $(M,\alpha,\Psi_\alpha)$ and $(M',\alpha',\Psi_{\alpha'})$. Show that $\varphi$ is a strict contactomorphism.

\item Consider $G = \mathbb{T}^2$ and thus $\mathfrak{g}^* = \mathbb{R}^2$ and $\mathbb{Z}_G = \mathbb{Z}^2$. We define $$v_1 = \begin{bmatrix}
0 \\ -1
\end{bmatrix} \quad \textrm{and} \quad v_2 =  \begin{bmatrix}
1 \\ 1
\end{bmatrix}.$$
Show that the set $$ C = \{\eta \in \mathbb{R}^2 | \langle \eta,v_1\rangle \geq 0, \langle \eta,v_2\rangle \geq 0\}$$ is a good rational polyhedral cone.

\item \label{oef questionable} Consider a symplectic manifold $(M,\omega)$ endowed with a Hamiltonian torus action. Let $\Phi$ be a corresponding (symplectic) momentum map. Let $S$ be a hypersurface of contact type of $(M,\omega)$ (with respect to the Liouville vector field $Y$) which is invariant under the torus action. Consider the contact form $\alpha = \iota_Y\omega|_S$ on $S$. Show that $\Phi|_S$ corresponds to the $\alpha$-momentum map for the induced torus action on $(S,\alpha)$.

\end{enumerate}

\end{document}